%% file: PAPER.tex
\documentclass{siamltex}
\usepackage{times}
\usepackage{graphics,epsfig,graphicx,color}
\usepackage{amsmath}
\usepackage{mathrsfs,amsfonts}
\usepackage{fullpage}
\usepackage{url}
\usepackage[debug,
letterpaper=true,
colorlinks=true,
linkcolor=red,
filecolor=green,
citecolor=red,
pdfpagemode=None]
{hyperref}

\graphicspath{{./Figs/}{./}}

\newtheorem{alg}[theorem]{Algorithm}

\newtheorem{rem}[theorem]{Remark}

\input{mathmacros.tex}

\title{Robust nonlinear processing of active array data in inverse
  scattering via truncated reduced order models} 
\author{
Liliana Borcea\footnotemark[1]
\and
Vladimir Druskin\footnotemark[2]
\and
Alexander V. Mamonov\footnotemark[3]
\and
Mikhail Zaslavsky\footnotemark[4]
}

\begin{document}

\maketitle

\renewcommand{\thefootnote}{\fnsymbol{footnote}}

\footnotetext[1]{Department of Mathematics, University of Michigan,
  Ann Arbor, MI 48109-1043 (borcea@umich.edu)}
\footnotetext[2]{Formerly Schlumberger-Doll Research Center, 
currently Druskin Algorithms (vdruskin@gmail.com)}
\footnotetext[3]{Department of Mathematics, University of Houston,
  Houston, TX 77204-3008 (mamonov@math.uh.edu)}
 \footnotetext[4]{Schlumberger-Doll Research Center, 1 Hampshire St.,
  Cambridge, MA 02139-1578 (mzaslavsky@slb.com)}

\begin{abstract}
We introduce a novel algorithm for nonlinear processing of data
gathered by an active array of sensors which probes a medium with
pulses and measures the resulting waves. The algorithm is motivated by
the application of array imaging. We describe it for a generic
hyperbolic system that applies to acoustic, electromagnetic or elastic
waves in a scattering medium modeled by an unknown coefficient called
the reflectivity. The goal of imaging is to invert the nonlinear
mapping from the reflectivity to the array data.  Many existing
imaging methodologies ignore the nonlinearity i.e., operate under the
assumption that the Born (single scattering) approximation is
accurate. This leads to image artifacts when multiple scattering is significant.
Our algorithm seeks to transform the array data to those corresponding to the Born
approximation, so it can be used as a pre-processing step for any
linear inversion method. The nonlinear data transformation algorithm
is based on a reduced order model defined by a proxy wave propagator operator 
that has four important properties. 
First, it is data driven, meaning that it is constructed from the data alone, 
with no knowledge of the medium. Second, it can be factorized in two operators
that have an approximately affine dependence on the unknown reflectivity. 
This allows the computation of the Fr\'{e}chet derivative of the reflectivity to the 
data mapping which gives the Born approximation. Third, the algorithm involves 
regularization which balances numerical stability and data fitting with accuracy of 
the order of the standard deviation of additive data noise.
Fourth, the algebraic nature of the algorithm makes it applicable
to scalar (acoustic) and vectorial (elastic, electromagnetic) wave data without any 
specific modifications.
\end{abstract}



%=================================================
\section{Introduction}
This paper introduces a nonlinear data processing algorithm motivated
by an inverse scattering problem for linear hyperbolic systems of
equations modeling acoustic, electromagnetic or elastic waves in a
heterogeneous, isotropic, non-attenuating medium. Specifically, we are
interested in array imaging, where a collection (array) of nearby
sensors probes the medium with pulses and measures the resulting
waves. These measurements, called the array data ${\bf D}$, are
processed in imaging to obtain an estimate of the medium, a.k.a. an
image.

Array imaging is an important technology in ocean acoustics
\cite{collins1994inverse}, nondestructive evaluation and structural
health monitoring \cite{achenbach2000quantitative}, diagnostic
ultrasound \cite{szabo2004diagnostic}, radar imaging
\cite{curlander1991synthetic,cheney2009fundamentals}, seismic imaging
\cite{claerbout1985imaging,aki1980quantitative,biondi20063d} and so
on. At the basic level, it seeks to estimate the medium, modeled by
unknown coefficients in the hyperbolic system, by minimizing in the
least squares sense the differences between the measured data ${\bf
  D}$ and the synthetic data predicted by the model. The mapping
between the coefficients and ${\bf D}$ is  nonlinear and
non-convex. Thus, iterative model updating using Newton-type
optimization methods is computationally demanding and stagnates at
local minima. Ever increasing computing power has brought progress
toward the solution of this nonlinear least squares problem, mostly in
the seismic community, where it is known as full-waveform inversion
\cite{virieux2009overview,brossier2015velocity,abubakar2011three,huang2017full}.
However, in many applications images need to be formed in real time,
so computationally intensive full-waveform inversion approaches cannot
be used. Furthermore, the lack of convexity of the objective function
remains a challenge, which is somewhat mitigated by
good initial guesses or processing data in expanded frequency bands,
from low to high \cite{bunks1995multiscale,virieux2009overview} see
also \cite{chen1997inverse}.

Due to these challenges, imaging remains largely based on a
combination of high frequency and linearization (Born) approximations,
where the medium is modeled as the sum of a smooth background and a
rough perturbation.  The wave speed $c$ of the background controls the
propagation of the waves through the medium, which is often modeled
with geometrical optics, whereas the rough perturbation, a.k.a. the
reflectivity, causes reflections which are registered at the array. In
applications like airborne imaging radar or non-destructive evaluation
of materials $c$ is known and constant. In seismic imaging $c$ is not
known and it is challenging to estimate it due to the oscillatory
nature of the waves, on the small scale of the wavelength, which
causes small perturbations of $c$ to result in large changes of the
waves.  This issue, known in the geophysics literature as cycle skipping \cite{virieux2009overview}, 
is at the heart of the lack of convexity of the least
squares data fitting functionals and the consequent stagnation of
iterative velocity updates at local minima that are physically
meaningless.  Specialized methodologies for determining $c$ have been
developed \cite{stolk2002smooth,uhlmann2001travel}, but they are
usually carried out separately from the estimation of the reflectivity
and require more data, gathered at large arrays.

The reflectivity estimation is commonly based on the Born
approximation, which assumes that the mapping between the rough part
of the medium properties and ${\bf D}$ is linear.  The linearization
of this mapping is studied in
\cite{beylkin1985imaging,BeylkinBurridge} and it is used in the
popular reverse time migration method
\cite{bleistein2013mathematics,biondi20063d} and the related filtered
back-projection \cite{cheney2009fundamentals} and matched filtering
\cite{therrien1992discrete,collins1994inverse} imaging
approaches. However, nonlinear (multiple scattering) effects are
present in ${\bf D}$ and these methods may produce images with
artifacts. The algorithm introduced in this paper seeks to transform
the array data gathered in strongly scattering media to those
corresponding to the Born approximation.  Thus, it can be used as a
nonlinear data pre-processing in conjunction with any linear imaging
algorithm.

The first question that arises when pursuing such a nonlinear
transformation is how to parametrize the medium i.e., how to define
the reflectivity function with respect to which we linearize. We base
our parametrization on the analysis in \cite{BeylkinBurridge} which
identifies the combinations of the medium parameters that give the
leading order contribution to the single scattered waves. This
contribution depends on the acquisition geometry, specifically on the angles
between the direction (rays) of the incoming and outgoing waves at the
array. We consider an array of small aperture size with respect to the
depth of the reflectors, so that these angles are small and the
leading contribution to the Born approximation of the data
is due to the variations of the logarithm of the wave impedance
$\sigma$ \cite[Sections 1,2]{BeylkinBurridge}. This is the unknown
reflectivity in our study, denoted by $q = \ln \sigma$, and we seek to
linearize the mapping $q \to {\bf D}$ under the assumption
that the wave speed $c$ is known.  In applications it is only the
smooth part of $c$ that is known or may be estimated via separate
velocity estimation \cite{stolk2002smooth,uhlmann2001travel}.
Nevertheless, we expect from the results in \cite{BeylkinBurridge}
that the unknown rough variations of $c$ alone will have a small
contribution in the Born approximation for small arrays.

In sonar array imaging the medium is described by a single wave speed
$c$ and impedance $\sigma$, defined in terms of the bulk modulus and
mass density. Similarly, the electric and magnetic wave fields
propagate at the same speed and there is a single impedance defined by
the electric permittivity and magnetic permeability. Thus, in both
acoustics and electromagnetics we have a single reflectivity $q$. In
elasticity the pressure (P) and shear (S) waves propagate at different
speeds $c_P$ and $c_S$ and there are two impedances $\sigma_P$ and
$\sigma_S$. These four functions are defined in terms of the mass
density and two Lam\'{e} parameters, so once we fix $c_P$ and $c_S$,
the two impedances depend on each other. This is why we still have a
single reflectivity function in our parametrization. The implication
is that we capture P-to-P and S-to-S scattering events in our
transformation, but we cannot resolve the P-to-S and S-to-P mode
conversions. This is consistent with the results in \cite[Section
  2]{BeylkinBurridge} which say that these mode conversion effects
become important only for large angles between the incoming and
outgoing rays i.e., for large arrays.

The idea of removing nonlinear, multiple-scattering effects from array
data has been pursued before. The studies in
\cite{borcea2010filtering,borcea2017time,borcea2011adaptive,alonso2011detection,
  aubry2009detection} propose various data filtering approaches for
improved imaging of point-like targets in random media. Filters of
multiply scattered waves in layered media are developed in
\cite{borcea2012filtering,fomel2007poststack}. The inversion of the
nonlinear reflectivity to data mapping using inverse Born and Bremmer
series methods has been proposed in
\cite{weglein1999multiple,malcolm2005method} for reflection seismology
and in \cite{moskow2012inverse} for optics. Boundary control methods
\cite{0266-5611-23-5-R01} and redatuming procedures
\cite{de2018exact,wapenaar2014marchenko,doi:10.1190/geo2015-0377.1}
have also been introduced.

Our algorithm is based on the data-driven reduced order model
approaches in \cite{druskin2016direct,druskin2018nonlinear,DtB}, which
are inspired by Krein's theory of Stieltjes strings
\cite{krein1951solution,kac1974spectral}. This theory has led to the
development of spectrally matched (optimal) grids that give spectrally
accurate finite difference approximations of Dirichlet-to-Neumann maps
\cite{druskin1999gaussian}. These grids have been used in inverse
problems in
\cite{borcea2011resistor,borcea2014model,borcea2005continuum} and play
a role in the reduced order modeling approach used in this paper and
in \cite{druskin2016direct,druskin2018nonlinear,DtB}. Related methods,
based on the theory of Marchenko, Gel'fand and Levitan as well as
Krein's
\cite{krein1951solution,krein1953transfer,marchenko1950,gelfand} have
been proposed for inverse scattering problems in layered media in
\cite{gopinath1971inversion,habashy1991generalized,burridge1980gelfand,
  symes1979inverse,santosa1982numerical,bube1983one} and have been
expanded recently to higher dimensions in
\cite{kabanikhin2011numerical,wapenaar2013three}. In this theory, the
inverse problem is reformulated in terms of nonlinear Volterra
integral equations. In the linear algebra setting this translates to
the Lanczos and Cholesky algorithms used in this paper and in
\cite{druskin2016direct,druskin2018nonlinear} or, alternatively, the
Stieltjes moment problems
\cite{gallivan1996some,gallivan1996rational,dyukarev2004indeterminacy}.

The construction of the reduced order model in this paper {follows
%and in
\cite{druskin2016direct,druskin2018nonlinear,DtB} and} is purely algebraic.
The reduced order model is defined by a proxy wave propagator that
maps the wave field from a given state at time $t$ to the future state
at time $t + \tau$, where $\tau$ is the time sampling interval of the
measurements at the array. The reduced model wave propagator has two
important properties: First, it is constructed from the array data,
using linear algebraic operations, without any knowledge of the
medium.  Second, it can be factored in two operators which are adjoint
to each other and have an affine dependence on the reflectivity
$q$. This allows us to calculate the Fr\'{e}chet derivative of the
mapping from the reflectivity to the data, and thus obtain the
linearized Born approximation map.

The algorithm in this paper is a robust version of the algorithm in
\cite{DtB}, which was developed in the context of imaging with sound
waves. The issue with the algorithm in \cite{DtB} is that it suffers
from numerical instability that can be controlled somewhat in one
dimension by a careful choice of the time sampling step
$\tau$. However, in higher dimensions the stability is also affected
by the sensor separation in the array. There is a trade-off between
spatial undersampling, which causes aliasing errors and oversampling
which introduces evanescent modes that cause instability. A good
sensor separation is at about half a wavelength, but in elasticity
different components of the waves propagate at different speed so a
good sampling for {shear waves corresponds to oversampling for pressure}
waves. This leads to ill-conditioned calculations which in combination
with unavoidable noisy data, cause the break-down of the algorithm in
\cite{DtB}. The regularization of the algorithm, using a reduced order
model constructed via spectral truncation, is the first main result of
the paper. The second result consists in showing how the algorithm can
be used {as a black-box} in array imaging with the three types of waves: sound,
electromagnetic and elastic.

The paper is organized as follows: In section \ref{sect:Alg} we begin
with the formulation of the problem, followed by a brief review of the
algorithm in \cite{DtB}. Then we introduce our robust data processing
algorithm. The presentation in section \ref{sect:Alg} is for a generic
linear hyperbolic system of equations written in symmetric form, in
terms of two first order partial differential operators that are
adjoint to each other and have an affine dependence on the unknown
reflectivity $q$. We also consider a symmetric array data model, and
write explicitly the mapping $q \to {\bf D}$. Sections
\ref{sect:Son}--\ref{sect:EL} are concerned with the application of
the algorithm to array imaging with sound, electromagnetic and elastic
waves. Specifically, they show how to derive the generic hyperbolic
system and data model used in Section \ref{sect:Alg}, starting from
the basic wave equation models: acoustic equations, Maxwell's
equations and elastic wave equations. In Section \ref{sect:Num} we
present numerical simulation results. We {conclude} with a summary in Section \ref{sect:sum}.

%=================================================
\section{Robust Data to Born transformation}
\label{sect:Alg}
In this section we give the robust algorithm for the nonlinear
transformation of the array data. We begin in section \ref{sect:form}
with the definition of the transformation, for a generic wave equation
satisfied by a vector valued wave field, in a setting that applies to
sound, electromagnetic and elastic waves.  To approximate numerically
the transformation, we use the data driven reduced order model defined
in section \ref{sect:revROM}. We review in section \ref{sect:rev.2}
the algorithm introduced in \cite{DtB} for approximating the
transformation and then describe its shortcomings. These motivate the
robust algorithm introduced in section \ref{sect:result}.

%=================================================
\subsection{Formulation of the problem}
\label{sect:form}
Consider the generic wave equation 
\begin{equation}
\partial^2_t  \bP^{(s)}(t,\bx) + L_q
  L_q^T\bP^{(s)}(t,\bx) = 0, \quad \bx \in \Omega, \quad t > 0,
\label{eq:1}
\end{equation}
for the vector valued wave field $\bP^{(s)}(t,\bx) \in \RR^{\dP}$,
where $\dP \ge 1$. This field is defined at time $t \ge 0$, at
locations $\bx$ in the half space domain\footnote{Other domains, such
  as the whole space, may be considered as well.}  $\Omega$ with
boundary $\partial \Omega$,
\[\Omega = \{ \bx = (x_j)_{1 \le j \le d}: ~
x_d > 0\}\subset \RR^d, \qquad \partial \Omega = \{ \bx = (x_j)_{1 \le
  j \le d}: ~ x_d = 0\},\] where $d \ge 1$.  The boundary is modeled
by some homogeneous boundary conditions satisfied by
$\bP^{(s)}(t,\bx)$.

The operator $L_q L_q^T$ in \eqref{eq:1} is symmetric, positive
definite. Its factors $L_q$ and $L_q^T$ are first order partial
differential operators with respect to $\bx$, adjoint to each other in
the Euclidian inner product. The index of these operators stands for
the reflectivity function $q(\bx)$, the unknown in the inverse
problem. This appears as a coefficient\footnote{Specific dependencies
of $L_q$ on $q$ are considered in sections \ref{sect:Son}--\ref{sect:EL}.} 
in the expressions of $L_q$ and $L_q^T$, which are affine in $q(\bx)$.

The waves are generated and measured by sensors at locations $\bx_l$,
for $l = 1, \ldots, m_{a}$, in a compact set of diameter $a$ on the
surface $ \{\bx = (x_j)_{1 \le j \le d}: x_d = 0+\} $ lying just above
$\partial \Omega$. The sensors are closely spaced so they behave like
a collective entity, called an active array of aperture size $a$.  The
index $(s)$ stands for the wave generated by the $s^{\rm th}$ source
excitation, modeled with the initial conditions
\begin{equation}
\bP^{(s)}(0,\bx) = \bb^{(s)}(\bx), \quad \partial_t \bP^{(s)}(0,\bx) =
0, \quad \bx \in \Omega.
\label{eq:2}
\end{equation}
In the case of scalar (sound) waves, $s$ indexes the location of the
sensors and it takes the values $s = 1, \ldots, m_{a}$. For vectorial
waves, $s$ also accounts for different polarizations of the
excitations. Thus, we let $s = 1, \ldots, m$ with $m \ge m_{a}$. In
either case, the excitation is modeled by the initial wave field
$\bb^{(s)}(\bx)$, which is compactly supported in the vicinity of the
source.

The resulting wave is measured at all the sensors, at discrete times
$t_k = k \tau$, with $k = 0, 1, \ldots, 2n -1$ and time increment
$\tau$.  For each $t_k$ we have an $m \times m$ data matrix\footnote{
  We explain in sections \ref{sect:Son}--\ref{sect:EL} how typical
  active array measurements can be put in the form \eqref{eq:3}.} with
entries indexed by $(r,s)$ and  modeled by
\begin{equation}
D^{(r,s)}_k = \int_{\Omega}{\rm d} \bx \, \bb^{(r)T}(\bx) \bP^{(s)}(k
\tau, \bx), \quad r,s = 1, \ldots, m.
\label{eq:3}
\end{equation}
Here $T$ denotes transpose, $r$ stands for the $r^{\rm th}$
measurement and $s$ for the $s^{\rm th}$ source excitation.

We gather all the  waves in the matrix valued field 
\begin{equation}
\bP(t,\bx) = \Big( \bP^{(1)}(t,\bx), \ldots, \bP^{(m)}(t,\bx) \Big)
\in \RR^{\dP \times m}.
\label{eq:4}
\end{equation}
This satisfies the wave equation
\begin{equation}
\partial_t^2 \bP(t,\bx) + L_q L_q^T \bP(t,\bx) =0 , \quad \bx \in
\Omega, \quad t > 0,
\label{eq:6}
\end{equation}
with some homogeneous boundary conditions at $\partial \Omega$ and the
initial conditions
\begin{equation}
\bP(0,\bx) = \bb(\bx), \quad \partial_t \bP(0,\bx) = 0,
\label{eq:7}
\end{equation}
defined by  the matrix 
\begin{equation}
\bb(x) = \Big( \bb^{(1)}(\bx), \ldots, \bb^{(m)}(\bx) \Big) \in
\RR^{\dP \times m}.
\label{eq:8}
\end{equation}
The differential operators in equation \eqref{eq:6} are understood to
act on one column of $\bP(t,\bx)$ at a time.

Let us introduce the wave propagator operator 
\begin{equation}
\cPq = \cos \Big( \tau \sqrt{L_q L_q^T} \Big),
\label{eq:dtb17}
\end{equation} 
which maps the wave from its initial state { at $t=0$} to the state at
time $\tau$,
\[
\bP(\tau,\bx) = \cos\Big(\tau \sqrt{L_q L_q^T}\Big) \bb(\bx) = \cPq
\,\bb(\bx).
\]
The wave at the $k^{\rm th}$ measurement instant $t_k = k \tau$,
called henceforth the $k^{\rm th}$ wave snapshot $\bP_k(\bx)$, is
given by
\begin{equation}
\bP_k(\bx) = \bP(k \tau, \bx) = \cos \Big( k \tau \sqrt{L_q L_q^T}
\Big)\bb(\bx),
\label{eq:8p}
\end{equation}
or equivalently, in terms of the propagator, by
\begin{equation}
\bP_k(\bx) = \cos \Big(k \arccos \cPq\Big) \bb(\bx) =
\cT_k(\cPq)\bb(\bx),
\label{eq:dtb18}
\end{equation}
where $\cT_k$ are the Chebyshev polynomials of the first kind. The
recurrence relations satisfied by these polynomials prove useful in
the construction of the reduced order model described in section
\ref{sect:revROM}.

The data matrices \eqref{eq:3} {written in terms of the propagator are}
%are
\begin{equation}
\bD_k = \int_{\RR^d}{\rm d} \bx \, \bb^T(\bx) \bP_k (\bx) = 
\int_{\RR^d}{\rm d} \bx \, \bb^T(\bx) \cT_k(\cPq)\bb(\bx), \quad
k = 0, 1, \ldots, 2n-1.
\label{eq:9}
\end{equation}
Although $L_q$ and $L_q^T$ are affine in the reflectivity $q$, it is
clear from \eqref{eq:dtb17} and \eqref{eq:9} that the mapping
\begin{equation}
q \to (\bD_k)_{0 \le k \le 2n-1}
\label{eq:map}
\end{equation} 
is nonlinear.  Most imaging algorithms are based on the assumption
that  \eqref{eq:map} may be approximated by the
linear map
\begin{equation}
q \to (\bD_k^{\rm Born})_{0 \le k \le 2n-1},
\label{eq:mapB}
\end{equation}
defined by the Fr\'{e}chet derivative of \eqref{eq:map} at $q = 0$.
This is the Born approximation, which is justified when the
reflectivity is small in some norm. The purpose of our algorithm is to
transform the data \eqref{eq:9}, acquired for a reflectivity $q$ that
is not small, to those corresponding to the linear Born model
\eqref{eq:mapB}.

We call henceforth the nonlinear transformation
\begin{equation}
(\bD_k)_{0 \le k \le 2n-1} \to \big(\bD_k^{\rm Born}\big)_{0 \le k \le 2n-1}
\label{eq:DtBmap}
\end{equation}
the Data to Born (DtB) mapping.

%=================================================
\subsection{The reduced order model}
\label{sect:revROM}

{Following \cite{DtB} we build our approach on the theory of data driven model order 
reduction. We refer to} the reduced model as data driven, because its construction
is based on matching the data \eqref{eq:9}. We define
it using the proxy wave propagator matrix $\cbPn \in \RR^{nm  \times nm}$ 
and the proxy {initial condition} matrix $\bbn \in \RR^{nm \times m}$,
satisfying the analogue of \eqref{eq:9},
\begin{equation}
\bD_k = \bb^{(n)T} \cT_k(\cbPn) \bb^{(n)}, \qquad k = 0, 1, \ldots, 2n-1.
\label{eq:ROM1}
\end{equation}

The data interpolation relations \eqref{eq:ROM1} are exact. However,
it is easier to explain the construction of $\cbPn$ and $\bbn$ if we
consider an approximation of the wave field \eqref{eq:4} by the matrix
$ \bP^{(N)}(t) \in \RR^{N \dP \times m}$ obtained from discretizing in
$\bx$ on a very fine grid with $N$ points\footnote{{
The fine grid discretization of $\bx$ is performed for expository reasons only 
and it is not required for the reduced model derivation. We refer to \cite{druskin2018nonlinear}
for the discretization-free derivation in the continuum.}}. 
Although $\bx$ lies in the half space $\Omega$, since the waves travel at 
finite speed, we can always restrict $\bx$ to a compact set 
$\Omega_c \subset \Omega$. We let $\Omega_c$ be a cube with side larger 
than the travel distance over the duration $2n \tau$, and with boundary defined 
by the union of two sets: The accessible boundary, lying in $\partial \Omega$ and 
the inaccessible boundary inside $\Omega$.  This extension allows us to impose
homogeneous Dirichlet boundary conditions at the inaccessible
boundary, without affecting the wave at $t \in [0,2 n \tau).$ The fine
  grid with $N$ points is the discretization of $\Omega_c$ and we
  group the unknowns in $\dP \times m$ blocks, ordered by the index of
  discretization of the coordinate $x_d$, starting from $x_d =0$ at
  $\partial \Omega$.  The discretization of the operators $L_q$ and
  $L_q^T$ gives the block lower bidiagonal matrix $\bL_q^{(N)}$ and
  its transpose $\bL_q^{(N)T}$, and the propagator operator
  \eqref{eq:dtb17} is approximated by the matrix
\begin{equation}
\cbPN = \cos \Big(\tau \sqrt{\bL_q^{(N)} \bL_q^{(N)T}}\, \Big) \in \RR^{N
  \dP \times N \dP}.
\label{eq:PropN}
\end{equation}
The initial field \eqref{eq:8} is
approximated by the matrix $\bb^{(N)} \in \RR^{N \dP \times m}$ and the
wave snapshots \eqref{eq:dtb18} become the matrices 
\begin{equation}
\bP_k^{(N)} = \cT_k(\cbPN) \bbN \in \RR^{N
  \dP \times m},
\label{eq:SnapN}
\end{equation}
for time indexes $k =0, 1, \ldots 2n-1$.

%=================================================
\subsubsection{The reduced order model as a projection}
\label{sect:ROMProj}
It is shown in \cite[Lemma 4.1]{druskin2018nonlinear} that $\cbPn$ and
$\bbn$ satisfying the data interpolation conditions \eqref{eq:ROM1}
can be defined as projections of the true propagator and initial condition matrices
\begin{equation}
\cbPn = \bV^T \cbPN \bV, \qquad \bbn = \bV^T \bbN,
\label{eq:ROM}
\end{equation}
{ on the block Krylov subspace 
\begin{equation}
\mbox{span}\big\{ \bbN, \cbPN \bbN, \ldots, \big(\cbPN\big)^{n-1} \bbN
\big\} = \mbox{range}\big( \bP^{(N)} \big),
\label{eq:ROM2}
\end{equation}
where $\bP^{(N)}$ is the matrix of the {first $n$} snapshots
\begin{equation}
\bP^{(N)} = \Big(\bP_0^{(N)}, \bP_1^{(N)}, \ldots \bP_{n-1}^{(N)}\Big)
\in \RR^{N \dP \times n m},
\label{eq:ROM3}
\end{equation}
and $\bV= (\bV_1, \ldots, \bV_n) \in \RR^{N \dP \times nm} $ is an
orthogonal matrix with columns spanning the subpsace \eqref{eq:ROM2}.
There are many such orthogonal matrices, but for our purpose we use
the one corresponding to the block QR factorization
\begin{equation}
\bP^{(N)} = \bV \bR, 
\label{eq:ROM4}
\end{equation}
which is shown in \cite[Sections 2, 3]{DtB} to be approximately
independent of the unknown reflectivity $q(\bx)$.  In this
factorization $\bV$ contains the orthonormal basis for the subspace
\eqref{eq:ROM2}, in which $\bP^{(N)}$ is transformed to the $n m
\times nm $ block upper tridiagonal matrix
\begin{equation}
\bR = \begin{pmatrix}
\bR_{1,1} & \bR_{1,2} & \ldots &\bR_{1,n} \\
0 & \bR_{2,2} & \ldots &\bR_{2,n} \\
\vdots\\
0 & 0 &\ldots  & \bR_{n,n} 
\end{pmatrix}, \qquad \bR_{i,j} \in \RR^{m \times m}, \quad 1 \le i \le j \le n.
\label{eq:ROM6}
\end{equation}

\vspace{0.05in}
\begin{rem}
\label{rem:2} The matrix $\bP^{(N)}$ is
already in almost block upper triangular form, because the wave 
equation is causal and the speed of propagation is finite. Explicitly, the first component
$\bP_0^{(N)}$, which is the initial condition $\bbN$, is supported in
the first blocks, corresponding to the points $\bx$ with component
$x_d \approx 0$. The second component $\bP_1^{(N)}$ is the wave at
time $\tau$, which advances a bit further in $x_d$, so it has a few
more nonzero blocks, and so on. This nearly block upper triangular
structure of $\bP^{(N)}$ is important because it gives that $\bV$,
which contains the orthonormal basis that transforms $\bP^{(N)}$ to
the block upper triangular $\bR$, is an approximate identity and is,
therefore, almost independent of the reflectivity $q(\bx)$.
\end{rem}

\vspace{0.05in}
\begin{rem}
\label{rem:1}
The components of $\bV$,
called the orthogonalized snapshots, satisfy the causality relations
\[
\bV_j \in \mbox{span}\big\{\bP_0^{(N)}, \ldots, \bP_{j-1}^{(N)}
\big\}, \quad j = 1, \ldots, n,
\]
derived from 
\begin{equation}
\bV = \bP^{(N)} \bR^{-1},
\label{eq:ROM4p}
\end{equation} 
and the block upper triangular structure of $\bR^{-1}$.
We also have from definition \eqref{eq:ROM} and the recursion
relations satisfied by the Chebyshev polynomials that $\cbPn$ has
block tridiagonal structure (see \cite[Appendix A]{DtB}).
\end{rem}

%=================================================
\subsubsection{Calculation of the reduced order model}

%The reduced order model $\big(\cbPn, \bbn\big)$ cannot be calculated
%as in \eqref{eq:ROM}, because neither $\cbPN$ nor $\bP^{(N)}$ are
%known. However, we can obtain it directly from the data \eqref{eq:9},
%as we now explain.

{The projection formulas \eqref{eq:ROM} cannot be used directly
to compute the reduced order model, since neither $\cbPN$ nor $\bP^{(N)}$ 
are known. However, it is possible to obtain the reduced model
just from the data \eqref{eq:9}, as we explain below.}

The calculation is based on the following {multiplicative property} 
of Chebyshev polynomials
\begin{equation}
\cT_i\big(\cbPN\big) \cT_j\big(\cbPN\big) = \frac{1}{2} \Big[
  \cT_{i+j}\big(\cbPN\big) + \cT_{|i-j|}\big(\cbPN\big)\Big].
\label{eq:ROM7}
\end{equation}
This identity and equation \eqref{eq:SnapN} give that the Gramian 
\begin{equation}
\bM = \bP^{(N)T} \bP^{(N)} = (\bM_{i,j})_{1 \le i, j \le n},
\label{eq:ROM8} 
\end{equation}
is defined by the data
as
\begin{equation}
\bM_{i,j} = \frac{1}{2} \Big(\bD_{i+j-2} + \bD_{|i-j|}\Big) \in \RR^{m
  \times m}, \qquad i,j = 1,\ldots, n.
\label{eq:14}
\end{equation}
Moreover, from \eqref{eq:ROM4} and the orthogonality of $\bV$ we get that 
\begin{equation}
\bM = \bR^T \bR,
\label{eq:ROM9}
\end{equation}
so $\bR$ can be computed using a block Cholesky factorization of
$\bM$.  We also obtain from \eqref{eq:ROM7} that the matrix
\begin{equation}
\bS = \bP^{(N)T} \cbPN \bP^{(N)} = (\bS_{i,j})_{1 \le i,j \le n},
\label{eq:14p}
\end{equation} 
is defined by the data as  
\begin{equation}
\bS_{i,j} = \frac{1}{4} \Big(\bD_{i+j-1} + \bD_{|j-i+1|} +
\bD_{|j-i-1|} + \bD_{|j+i-3|} \Big) \in \RR^{m \times m}, \qquad i,j =
1,\ldots, n.
\label{eq:15}
\end{equation}

The reduced order model is obtained from \eqref{eq:ROM} and
\eqref{eq:ROM4p}, 
\begin{align}
\cbPn &= \bR^{-T} \bS \bR^{-1},
\label{eq:16}\\
  \bbn & = {
  \bR^{-T} \bP^{(N) T} \bbN = 
   \bR^{-T} \bP^{(N) T} \bP^{(N)} \bE_1 } =
  \bR \bE_1,
\label{eq:16p}
\end{align}
where we introduced the matrix
\begin{equation}
  \bE_1 = \begin{pmatrix} {\bf I}_m \\ {\bf 0}_m \\ \vdots \\ {\bf
      0}_m
  \end{pmatrix} \in \RR^{nm \times m},
  \label{eq:16pp}
\end{equation}
with ${\bf I}_m$ the $m \times m$ identity and ${\bf 0}_m$ the $m
\times m$ identically zero matrix. We also have from \eqref{eq:ROM1}
and \eqref{eq:16p} that
\begin{equation}
\bD_0 = \bb^{(n)T} \bb^{(n)} = \bR_{1,1}^T \bR_{1,1}.
\label{eq:16ppp}
\end{equation}

%=================================================
\subsubsection{Factorization of the wave  propagator}
To explain the algorithm for approximating the DtB transformation
\eqref{eq:DtBmap}, we use the factorization
\begin{equation}
\frac{2}{\tau^2} ({\bf I}_{_{N \hspace{-0.01in}\dP}} - \cbPN) = \cbLN_q \cbLNT_q,
\label{eq:R1}
\end{equation}
where ${\bf I}_{_{N \hspace{-0.01in}\dP}}$ is the $N \dP \times N \dP$ identity matrix and
\begin{align}
  \cbLN_q &= \frac{2}{\tau} \bL^{(N)}_q (\bL_q^{(N)T}
  \bL_q^{(N)})^{-1/2} \sin \Big(\frac{\tau}{2} \sqrt{\bL_q^{(N)T}
    \bL_q^{(N)}}\Big). \label{eq:R2}
\end{align}
This can be checked easily, using definition \eqref{eq:PropN} and the
singular value decomposition of $\bL^{(N)}_q$. By the definition of
$\cbPN$ in terms of the cosine, the left hand side in \eqref{eq:R1} is
positive semidefinite. We can make it positive definite for a small
time sample interval $\tau$, satisfying the Courant-Friedrich-Levy (CFL) type condition
\begin{equation}
{ \tau < \| \bL_q^{(N) T} \bL_q^{(N)} \|_{2}^{-1/2}, }
\label{eq:R3}
\end{equation}
and then obtain from the Taylor expansion of \eqref{eq:R2} that
\begin{equation}
  \cbLN_q \approx \bL_q^{(N)}.
\label{eq:R4}
\end{equation}
Thus, the matrix factor $\cbLN_q$ is an approximation of
$\bL_q^{(N)}$, which has lower block bidiagonal structure and has an
affine dependence on the unknown reflectivity $q$.

The reduced order model propagator has  a similar factorization 
\begin{equation}
\frac{2}{\tau^2} ({\bf I}_{_{nm}} - \cbPn) = \cbLn_q \cbLnT_q,
\label{eq:R5}
\end{equation}
with ${\bf I}_{_{nm}}$ the $nm \times nm$ identity matrix.  From
equations \eqref{eq:ROM} and \eqref{eq:R1} we have
\begin{equation}
\cbLn_q \cbLnT_q = \bV^T  \cbLN_q \cbLNT_q \bV.
\label{eq:R6}
\end{equation}
Furthermore, it is shown in \cite[Section 2.4]{DtB} that 
\begin{equation}
\cbLn_q = \bV^T  \cbLN_q \widehat{\bV},
\label{eq:R7}
\end{equation}
where $\widehat{\bV}$ is the orthogonal matrix {containing the orthonormal basis for}
the subspace spanned by the time snapshots of the dual wave
$\widehat{\bP}^{(N)}(t)$ at $t = (k+1/2) \tau$, for $k = 0, \ldots,
2n-1$.  This wave appears in the first order system formulation of the
wave equation
\begin{equation}
\partial_t \begin{pmatrix} 
\bP^{(N)}(t) \\
\widehat\bP^{(N)}(t)
\end{pmatrix} = \begin{pmatrix} {\bf 0} & - \bL_q^{(N)} \\
\bL_q^{(N)T} & {\bf 0} \end{pmatrix}\begin{pmatrix} 
\bP^{(N)}(t) \\
\widehat\bP^{(N)}(t)
\end{pmatrix},
\label{eq:R8}
\end{equation}
with initial conditions 
\begin{equation}
\bP^{(N)}(0) = \bbN, \quad \widehat{\bP}^{(N)}(0) = {\bf 0},
\label{eq:R9}
\end{equation}
and the construction of $\widehat{\bV}$ is similar to that of $\bV$.

\vspace{0.05in}
\begin{rem}
  \label{rem:3}
The reduced order model matrix \eqref{eq:R7} corresponds to the {Galerkin} approximation of the system
\eqref{eq:R8}--\eqref{eq:R9}, on the spaces of the primary and dual
snapshots, using the orthogonal bases in $\bV$ and $\widehat{\bV}$.
Recall from Remark \ref{rem:1} that $\cbPn$ has block tridiagonal
structure. Therefore, its Cholesky factor $\cbLn_q$ is block lower
bidiagonal. By equation \eqref{eq:R4} $\cbLN_q$ is approximately block
lower bidiagonal. The projection in \eqref{eq:R7}, with the orthogonal
matrices $\bV$ and $\widehat \bV$, maps $\cbLN_q$ to the lower block
bidiagonal $\cbLn_q$.
\end{rem}

%=================================================
\subsubsection{The reduced order model snapshots}

The reduced order model propagator $\cbPn$ and initial state $\bbn$
define the reduced order model snapshots,
\begin{equation}
\bP_k^{(n)} = \cT_k(\cbPn)\bbn,
\label{eq:Snapn}
\end{equation}
the analogues of \eqref{eq:SnapN}. These satisfy the algebraic system
of equations
\begin{align}
&
\frac{\bP^{(n)}_{k+1} - 2 \bP_k^{(n)} + \bP_{k-1}^{(n)}}{\tau^2} + \cbLn_q \cbLnT_q
{ \bP^{(n)}_k }= 0, 
\quad k \ge 0,
\nonumber \\
&\bP_0^{(n)}= \bbn = \bR \bE_1, \qquad \bP^{(n)}_{-1} = \bP_1^{(n)},
\label{eq:dtb22}
\end{align}
derived from \eqref{eq:R5} and the recursion relation of the Chebyshev
polynomials
\begin{equation}
\cT_{k+1}(\cbPn) + \cT_{k-1}(\cbPn) = 2 \cbPn \cT_k(\cbPn)~~
\mbox{for}~ k \ge 1, \label{eq:dtb19} \qquad \cT_0(\cbPn) = {\bf
  I}_{mn}, \quad \cT_1(\cbPn) = \cbPn.
\end{equation}

\vspace{0.05in}
\begin{rem}
\label{rem:FD1}
The discrete system \eqref{eq:dtb22} is designed to match exactly the
data
\begin{equation}
\bD_k = \bb^{(n)T} \bP_k^{(n)}, \quad k = 0, 1, \ldots, 2n-1,
\label{eq:R16}
\end{equation}
as follows from \eqref{eq:ROM1} and \eqref{eq:Snapn}. Moreover, from
\eqref{eq:R4}, \eqref{eq:R7} and Remarks \ref{rem:2}, \ref{rem:3} we
conclude that $\cbLn_q$ inherits approximately the affine dependence
of $\bL_q^{(N)}$ on the unknown reflectivity $q$.
\end{rem}

\vspace{0.05in}
\begin{rem}
\label{rem:FD}
The system \eqref{eq:dtb22} is a finite difference discretization of
the wave equation \eqref{eq:6}, with the time derivative
$\partial_t^2$ replaced by the central difference with time step
$\tau$ and the first order partial differential operator $L_q$
replaced by the lower block bidiagonal matrix $\cbLn_q$. 
The block structure of $\cbLn_q$ corresponds to a two point finite
difference discretization of the derivative $\partial_{x_d}$, whereas
the block Cholesky calculation of $\cbLn_q$ with \cite[Algorithm
  4.1]{DtB} corresponds to the discretization of the 
{remaining components} of the gradient 
$(\partial_{x_j})_{1 \le j \le d-1}$. 

It is shown in \cite[Sections 3.1, 4.6]{DtB} that the resulting finite
difference scheme is approximately the same as in the reference medium
with reflectivity zero i.e., on almost the same grid. This is a 
special grid that is closely connected with the reduced order
model. Moreover, the unknown reflectivity $q$ may be approximated on
this grid from the entries in the matrix $\cbLn_q$.
\end{rem}

%=================================================
\subsection{Review of the  DtB algorithm}
\label{sect:rev.2}
We now review the algorithm introduced in \cite{DtB} for approximating
the DtB transformation \eqref{eq:DtBmap}.  We begin by rewriting the
data matching relations \eqref{eq:ROM1} using the block Cholesky
factorization \eqref{eq:R1}, with $\cbPn$ and $\bbn$ calculated from 
\eqref{eq:16}--\eqref{eq:16p},
\begin{equation}
\bD_k = \bb^{(n)T} \cT_k\Big({\bf I}_{mn} - \frac{\tau^2}{2} \cbLn_q
\cbLnT_q\Big)\bb^{(n)}, \quad k = 0, 1, \ldots, 2n-1.
\label{eq:R10}
\end{equation}
We already explained that $\cbLn_q$ is approximately affine in $q$ and
we conclude from Remark \ref{rem:2} and \eqref{eq:ROM} that $\bbn$ is
approximately independent of $q$.  The Born data model \eqref{eq:mapB}
is given by the Fr\'{e}chet derivative of the mapping $q \to \big(
\bD_k \big)_{0 \le k \le 2n-1}$ evaluated at $q = 0$.

To be more explicit, let 
\begin{equation}
\bD_{0,k} = \bb^{(n)T} \cT_k\Big({\bf I}_{mn} - \frac{\tau^2}{2} \cbLn_0
\cbLnT_0\Big)\bb^{(n)}, \quad k = 0, 1, \ldots, 2n-1,
\label{eq:R11}
\end{equation}
be the analogue of \eqref{eq:R10} in the reference medium with $q =
0$.  These are called the reference data although they are not
measured, but are calculated analytically or numerically. We index
them by $0$ to distinguish them from the measurements 
{corresponding to the medium with reflectivity $q$}.
The analogue of \eqref{eq:R5} in the reference medium is
\begin{equation}
\frac{2}{\tau^2} ({\bf I}_{_{nm}} - \cbPno) = \cbLn_0 \cbLnT_0,
\label{eq:R12}
\end{equation}
with the reduced order model propagator $\cbPno$ calculated using
equations \eqref{eq:ROM8}--\eqref{eq:16} and $(\bD_{0,k})_{0 \le k \le
  2n-1}.$ Let {$\ep q$} be the scaled unknown reflectivity,
with small and positive $\ep$. The operator {$\cbLn_{\ep q}$}
corresponding to this scaled reflectivity is approximately
\begin{equation}
{\cbLn_{\ep q}} \approx \cbLn_0 + \ep \Big(\cbLn_q-\cbLn_0\Big),
\label{eq:R13}
\end{equation}
because $\cbLn_q$ is approximately affine in $q$. Then, the Born data
can be approximated as
\begin{equation}
\bD_{k}^{\rm Born} \approx \bD_{0,k} + \bb^{(n)T} \frac{d}{d \ep}
\cT_k\Big({\bf I}_{mn} - \frac{\tau^2}{2} \big[\cbLn_0 + \ep
  \big(\cbLn_q-\cbLn_0\big)\big] \big[\cbLn_0 + \ep
  \big(\cbLn_q-\cbLn_0\big)\big]^T\Big)\Big|_{\ep = 0} \bb^{(n)},
\label{eq:R14}
\end{equation}
for $k = 0, 1, \ldots, 2n-1$.

The derivative in \eqref{eq:R14} cannot be obtained with the chain
rule, because matrices do not commute with their derivatives. The
calculation of the derivative is in \cite[Algorithm 2.2]{DtB}. Here we
write it more explicitly, so that we can extend it to the robust DtB
algorithm. We begin by rewriting \eqref{eq:dtb22} in the first order
system form:
\begin{align}
 \frac{\bP_{k+1}^{(n)}-\bP_k^{(n)}}{\tau} &= - \cbLn_q \widehat
 \bP_k^{(n)}, \quad k = 0, \ldots, 2n-2, \label{eq:R17}\\
 \frac{\widehat
   \bP_{k}^{(n)}-\widehat \bP_{k-1}^{(n)}}{\tau} &= 
\cbLnT_q\bP_k^{(n)}, \quad k = 1, \ldots,  2n-1, \label{eq:R18}
\end{align}
with initial conditions 
\begin{equation}
\bP_0^{(n)} = \bbn, \quad \widehat \bP_0^{(n)} + \widehat \bP_{-1}^{(n)} = {\bf 0}.
\label{eq:R19}
\end{equation}
This is the reduced order model discretization of the first order
system \eqref{eq:R8}--\eqref{eq:R9}, with the large matrix $\cbLN_q$
replaced by the much smaller lower block bidiagonal matrix $\cbLn_q$.

In the DtB transformation \eqref{eq:R14} the matrix {$\cbLn_{\ep q}$}
for the scaled reflectivity {$\ep q$} is approximated by
\eqref{eq:R13}, so the right hand sides in
\eqref{eq:R17}--\eqref{eq:R18} become affine in $\ep$. Thus, we obtain
from \eqref{eq:Snapn} and \eqref{eq:R16} that
\begin{equation}
\bb^{(n)T} \frac{d}{d \ep} \cT_k\Big({\bf I}_{mn} - \frac{\tau^2}{2}
\big[\cbLn_0 + \ep \big(\cbLn_q-\cbLn_0\big)\big] \big[\cbLn_0 + \ep
  \big(\cbLn_q-\cbLn_0\big)\big]^T\Big)\Big|_{\ep = 0} \bb^{(n)} =
\bb^{(n)T}\delta \bP_k^{(n)}, 
\label{eq:R21}
\end{equation}
for $k = 0, \ldots, 2n-1$, with {the snapshot perturbation} $
\delta \bP_k^{(n)}$ calculated from the time stepping scheme
\begin{align}
 \frac{\delta \bP_{k+1}^{(n)}-\delta \bP_k^{(n)}}{\tau} + \cbLn_0
 \delta \widehat \bP_k^{(n)} &= - \big(\cbLn_q - \cbLn_0 \big) \hat
 \bP_{0,k}^{(n)}, 
 \quad k = 0, \ldots, 2n-2, \label{eq:R22}\\ 
\frac{\delta \widehat \bP_{k}^{(n)} - 
 \delta \widehat \bP_{k-1}^{(n)}}{\tau} - \cbLnT_0 \delta  \bP_k^{(n)} &=
 \big(\cbLnT_q - \cbLnT_0 \big)  \bP_{0,k}^{(n)},
 \quad k = 1, \ldots, 2n-1, \label{eq:R23}
\end{align}
with homogeneous initial conditions 
\begin{equation}
\delta \bP_0^{(n)} = {\bf 0}, \quad \delta \widehat \bP_0^{(n)} +
\delta \widehat \bP_{-1}^{(n)} = {\bf 0}.
\label{eq:R24}
\end{equation}

%=================================================
\subsubsection{Shortcomings of the DtB algorithm}
\label{sect:rev.3}

The main issue with the calculation of \eqref{eq:R14} is that the
Gramian $\bM$ has poor condition number, due to a few very small
eigenvalues. The construction of the reduced order model is based on
the Cholesky factorization {\eqref{eq:ROM9} of the Gramian},
and it breaks down when $\bM$ obtained from \eqref{eq:ROM8} becomes 
indefinite due to noisy data.

For acoustic waves and in one dimension, it is shown in \cite[Section
  6]{druskin2016direct} that the condition number of $\bM$ depends on
the time sampling rate $\tau$. A good choice of $\tau$ corresponds to
the Nyquist sampling rate for the temporal frequencies of the pulse
emitted by the sources.  A much smaller time step gives a very poor
condition number, because the wave snapshots at two consecutive times
are nearly the same, and a larger time step is not desirable because
it causes temporal aliasing errors. It is also shown in
\cite{druskin2016direct} that not all $\tau$ similar to the Nyquist
sampling rate are the same, so $\tau$ should be selected using
numerical calibration. This impedes an automated data processing
procedure.

In multiple dimensions the condition number of $\bM$ also depends on
the spatial frequency (wavenumber) of the measurements, determined by
the spacing of the sensors in the array. From Fourier transform theory
we know that we can resolve waves with wavenumber $\kappa$ of at most
$\pi/h$, when the sensors are spaced at distance $h$ in the array
aperture. To minimize spatial aliasing errors in imaging, a typical
choice is $h = \lambda_c/2$, where $\lambda_c$ is the carrier (central)
wavelength of the source signals, so the waves are resolved up to the
wavenumber $\kappa = 2 \pi/\lambda_c$. 

To see what this means, suppose for a moment that the medium where
homogeneous, so we could decompose the wave field in independent plane
waves by Fourier transforming in $t$ and $(x_j)_{1 \le j \le
  d-1}$. These waves are of the form $\exp( i \boldsymbol{\kappa}^\pm
\cdot \bx - i \omega t)$, where $\omega$ is the temporal frequency and
$\boldsymbol{\kappa}^\pm = \big((\kappa_j)_{1 \le j \le d-1}, \pm
\kappa_d \big)$ is the wave vector with components
\begin{equation}
\label{eq:PW}
|\kappa_j| \le \frac{\pi}{h} ~ ~ \mbox{for} ~~ 1 \le j \le d-1, \quad
\kappa_d = \Big[\Big(\frac{2 \pi}{\lambda}\Big)^2 - \sum_{j=1}^{d-1}
  \kappa_j^2 \Big]^{1/2}.
\end{equation}
Here $\lambda$ is the wavelength at frequency $\omega$, which
typically satisfies $\lambda \approx \lambda_c$, and the $\pm$ sign
corresponds to forward and backward going waves along the direction
$x_d$.  In heterogeneous media such a decomposition still holds, but
the waves interact with each other due to scattering. 

We see from \eqref{eq:PW} that when the sensors are spaced at $h \le
\lambda/2$, we have evanescent modes, with imaginary $\kappa_d$, and
waves that propagate slowly in the $x_d$ direction, with real and small $\kappa_d$. These give
small contributions to the matrix \eqref{eq:ROM3}, corresponding
to the right singular vectors $\bzet_j$ of $\bP^{(N)}$ for nearly zero
singular values $\sigma_j$. The Gramian matrix \eqref{eq:ROM8} has the
eigenvectors $\bzet_j$ and eigenvalues $\sigma_j^2$, so it is poorly
conditioned.  One could try to control its condition number while
minimizing aliasing errors, by choosing $h$ a bit larger than
$\lambda_c/2$. However, in elasticity, different types of waves
propagate at different speed, so a good spatial sampling for the
{shear waves corresponds to oversampling for the pressure} waves.

All {the aforementioned} effects combined lead to unavoidable 
poor conditioning of $\bM$, which causes instability of the approximation 
\eqref{eq:R14} of the DtB mapping. The robust algorithm introduced in the next 
section regularizes this transformation using a spectral truncation.

%=================================================
\subsection{Robust DtB algorithm}
\label{sect:result}
There are two points to address in the regularization of the
approximation \eqref{eq:R14}: The first is the construction of the
regularized data map, via spectral truncation of the
Gramian $\bM$. This is explained in section \ref{sect:Trunc}. The
second is the calculation of the Fr\'{e}chet derivative of the
regularized data mapping, explained in section \ref{sect:Frech}. The
robust algorithm is summarized in section \ref{sect:SumAlg}.

%=================================================
\subsubsection{Regularized data map}
\label{sect:Trunc}
Let
\begin{equation}
\bM = {\itbf Z} \boldsymbol{\Sigma}^2 {\itbf Z}^T, \qquad {\itbf Z} =
(\bzet_1, \ldots, \bzet_{nm}), \quad \boldsymbol{\Sigma}^2 =
\mbox{diag} (\sigma_1^2, \ldots, \sigma_{nm}^2),
\label{eq:eigM}
\end{equation}
be the eigenvalue decomposition of the Gramian $\bM$, with eigenvalues
$\sigma_j^2$ listed in decreasing order.  To regularize the
calculation of the DtB transformation, we filter out its eigenvectors
$\bzet_j$ for eigenvalues $\sigma_j^2 \le \theta$, using the
orthogonal matrix
\begin{equation}
\bZ = (\bzet_1, \ldots, \bzet_{zm}) \in \RR^{n m \times zm}.
\label{eq:defZ}
\end{equation}
Here $\theta$ is of the order of the standard deviation of the noise
and $z$ is a natural number satisfying $z \le n$, chosen so that
$\sigma_j^2 > \theta$, for $j = 1, \ldots, zm$. We keep the
dimension of the projection space a multiple of $m$, because we wish
to use the block structure of the projected matrices, with blocks of
size $m \times m$.

To write the effect of the truncation on the system \eqref{eq:dtb22},
used in the calculation of the derivative of the regularized data
map,  we need two steps:

First, we rewrite \eqref{eq:dtb22} in terms of the matrices $\bM$ and
$\bS$ that can be computed directly from the data, using 
\begin{equation}
\bet_k = \bR^{-1} \bP_k^{(n)} \in \RR^{nm \times m}.
\label{eq:R25}
\end{equation}
Multiplying \eqref{eq:dtb22} on the left by $\bR^T$ we get from
definitions \eqref{eq:ROM9}--\eqref{eq:14p} that
\begin{align}
&\bM \frac{(\bet_{k+1}-2 \bet_k + \bet_{k-1})}{\tau^2} +
\frac{2}{\tau^2} (\bM - \bS) \bet_k = 0, \qquad k \ge 0, \nonumber \\
&\bet_0 = \bE_1, \quad \bet_{-1} = \bet_1.
\label{eq:R26}
\end{align}
We also obtain from \eqref{eq:16p} and \eqref{eq:R16} that the data
are related to \eqref{eq:R25} by
\begin{equation}
\bD_k = \bE_1^T \bM \bet_k, \qquad k = 0, \ldots, 2n-1.
\label{eq:R28}
\end{equation}

Second, because $\bM$ is likely indefinite due to noise, we define the
positive definite matrix
\begin{equation}
\bSig^2 = \bZ^T \bM \bZ = \mbox{diag}(\sigma_1^2, \ldots, \sigma_{zm}^2),
\label{eq:R29}
\end{equation}
and the matrix
\begin{equation}
\cbS = \bZ^T \bS \bZ \in \RR^{zm \times zm}.
\label{eq:R30}
\end{equation}
The solution $\bett_k$ of
\begin{align}
&\bSig^2 \frac{(\bett_{k+1}-2 \bett_k + \bett_{k-1})}{\tau^2} +
  \frac{2}{\tau^2} (\bSig^2 - \cbS) \bett_k = 0, \qquad k \ge
  0, \nonumber \\ &\bett_0 = \widetilde \bE_1 = \bZ^T \bE_1, \quad
  \bett_{-1} = \bett_1.
\label{eq:R34}
\end{align}
is approximately 
\begin{equation}
\bett_k \approx \bZ^T \bet_k,
\label{eq:R26pp}
\end{equation}
as can be seen by replacing $\bett_k$ by $\bZ^T \bet_k$ in the left
hand side of \eqref{eq:R34}, and then using equations \eqref{eq:R26},
\eqref{eq:R29} and the approximation $\bZ \bZ^T \bet_k \approx
\bet_k$, with error depending on the truncation
threshold\footnote{A detailed description of reduced order models
  based on spectral truncations is in
  \cite{antoulas2001survey,gugercin2004survey}.}  $\theta$.

As in any regularization scheme, we have a trade-off between the
stability and data fit. Instead of \eqref{eq:R28}, we have
\begin{equation}
\widetilde \bD_k = \widetilde \bE_1^T \bSig^2 \bett_k = \bE_1^T \bZ
\bZ^T \bM \bZ \bett_k \approx \bE_1^T \bZ \bZ^T \bM \bZ \bZ^T \bet_k,
\label{eq:R32}
\end{equation}
where the norm of the misfit $\bD_k - \widetilde \bD_k$ is bounded in
terms of the sum of the smallest $(n-z)m$ eigenvalues $\sigma_j^2$ of
$\bM$ \cite{antoulas2001survey}. Since $\sigma_j^2 \le \theta$, the
data are fit with an accuracy commensurate with the standard deviation
of the noise. 

We call 
\begin{equation}
q \to \Big(\widetilde \bD_k\Big)_{0 \le k \le 2n-1}
\label{eq:R32p}
\end{equation}
the regularized data map.

%=================================================
\subsubsection{Derivative of the regularized data map}
\label{sect:Frech}

To calculate the derivative of \eqref{eq:R32p}, we proceed as in
section \ref{sect:rev.2} and obtain from equation \eqref{eq:R34} a
first order system, the analogue of \eqref{eq:R17}--\eqref{eq:R18}.

Define 
\begin{equation}
\bpi_k^{(n)} = \bSig \bett_k \in \RR^{z m \times m},
\label{eq:R37}
\end{equation}
and multiply \eqref{eq:R34} on the left by $\bSig^{-1}$ to obtain the
system
\begin{align} 
&\frac{\bpi_{k+1}^{(n)} - 2 \bpi_k^{(n)} + \bpi_{k-1}^{(n)}}{\tau^2} +
  \frac{2}{\tau^2} \Big({\bf I}_{zm} - \bSig^{-1} \cbS \bSig^{-1}
  \Big) \bpi_k^{(n)} = 0, \qquad k \ge 0, \nonumber \\ &
  \bpi_0^{(n)} = \bSig \widetilde{\bE}_1, \quad \bpi_{-1}^{(n)} =
  \bpi_1^{(n)}.
\label{eq:R39} 
\end{align}
This is similar to \eqref{eq:dtb22}, except that instead of the
reduced order model propagator $\cbPn$ we have $ \bSig^{-1} \cbS
\bSig^{-1}$.  To compare these two matrices, recall equation
\eqref{eq:ROM9} and let
\begin{equation}
\bR = {\itbf U} \boldsymbol{\Sigma} {\itbf Z}^T, 
\label{eq:N1}
\end{equation}
be the singular value decomposition of the block upper triangular
matrix $\bR$, with ${\itbf Z}$ and $\boldsymbol{\Sigma}$ defined by
the eigenvalue decomposition \eqref{eq:eigM} of the Gramian, and with
orthogonal matrix ${\itbf U} \in \RR^{nm \times nm}$.  Substituting
\eqref{eq:N1} in \eqref{eq:16}, we obtain that
\begin{equation}
\cbPn = \bR^{-T} \bS \bR^{-1} = {\itbf U} \boldsymbol{\Sigma}^{-1}
      {\itbf Z}^T \bS {\itbf Z} \boldsymbol{\Sigma}^{-1} {\itbf U}^T,
\label{eq:N2}
\end{equation}
whereas
\begin{equation}
\bSig^{-1} \cbS \bSig^{-1} = \bSig^{-1} \bZ^T \bS \bZ \bSig^{-1}.
\label{eq:N3}
\end{equation}
These expressions are similar, except that the square $nm \times n m$
matrices $\boldsymbol{\Sigma}$ and ${\itbf Z}$ are replaced by the
truncations $\bSig \in \RR^{zm \times zm}$ and $\bZ \in \RR^{nm \times
  zm}$, and ${\itbf U}$ is replaced by the $zm \times zm$
identity. Because of these replacements, the matrix $\bSig^{-1} \cbS
\bSig^{-1}$ is not block tridiagonal. As we explained in Remark
\ref{rem:FD}, the block tridiagonal structure of $\cbPn$ is important,
because it corresponds to a finite difference scheme of the wave
equation \eqref{eq:1} with the discretization of the operator $L_q$
that is approximately affine in the unknown reflectivity $q$. 

{In order to recover the block tridiagonal structure we apply
the block Lanczos algorithm \cite{golub1977block,grimes1994shifted} to the
matrix $\bSig^{-1} \cbS \bSig^{-1}$ and the initial block 
$ \bpi_0^{(n)} = \bSig  \widetilde \bE_1 = \bSig \bZ^T \bE_1 \in \RR^{zm \times m}$.
It produces the orthogonal change of coordinates matrix 
$\boldsymbol{\mathcal U} \in \RR^{zm \times zm}$ such that}
%Thus, we change coordinates using the orthogonal matrix
%$\boldsymbol{\mathcal U} \in \RR^{zm \times zm}$, calculated with the
%block Lanczos algorithm \cite{golub1977block,grimes1994shifted} to
%transform \eqref{eq:N3} to the block tridiagonal matrix
\begin{equation}
\tilde{\boldsymbol{\mathscr P}}_q^{(n)} = \boldsymbol{\mathcal U}
\bSig^{-1} \cbS \bSig^{-1} \boldsymbol{\mathcal U}^T = 
{\boldsymbol{\mathcal U} \bSig^{-1} \bZ^T \bS 
\bZ \bSig^{-1} \boldsymbol{\mathcal U}^T},
\label{eq:N4}
\end{equation}
{is block tridiagonal. We call $\tilde{\boldsymbol{\mathscr P}}_q^{(n)}$ }
the regularized reduced order model propagator. The corresponding 
snapshots are 
\begin{equation}
\widetilde{\bP}^{(n)}_k = \boldsymbol{\mathcal U} \bpi_k^{(n)}, \qquad
k \ge 0,
\label{eq:N5}
\end{equation}
and they satisfy the algebraic system
\begin{align} 
& \frac{\widetilde{\bP}_{k+1}^{(n)} - 2\widetilde{\bP}_{k}^{(n)} +
    \widetilde{\bP}_{k-1}^{(n)}}{\tau^2} + \frac{2}{\tau^2} \Big( {\bf
    I}_{zm} - \tilde{\boldsymbol{\mathscr P}}_q^{(n)}\Big)
  \widetilde{\bP}_{k}^{(n)} = 0, \qquad k \ge 0, \nonumber \\ &
  \widetilde{\bP}^{(n)}_0 = \widetilde{\bb}^{(n)}, \quad
  \widetilde{\bP}_{-1}^{(n)} = \widetilde{\bP}_1^{(n)}, \label{eq:N6}
\end{align}
with initial condition
\begin{equation}
\widetilde{\bb}^{(n)} = \widetilde{\bR} \bE_1,
\label{eq:N6p}
\end{equation}
 Note that this is similar to
\begin{equation}
\bbn = \bR \bE_1 = {\itbf U} \boldsymbol{\Sigma} {\itbf Z}^T {\bf E}_1,
\label{eq:N6pp}
\end{equation}
except that $\bR$, with the singular value decomposition
\eqref{eq:N1}, is replaced by  (recall \eqref{eq:16p},
\eqref{eq:R34} and \eqref{eq:R39})
\begin{equation}
\widetilde{\bR} = \boldsymbol{\mathcal U}\bSig \bZ^T.
\label{eq:N6ppp}
\end{equation}
Note also that the recursion relations \eqref{eq:dtb19} of the
Chebyshev polynomials give that the regularized snapshots, the
solution of \eqref{eq:N6}, are 
\begin{equation}
\widetilde{\bP}_{k}^{(n)} = \cT_k\Big(\tilde{\boldsymbol{\mathscr
    P}}_q^{(n)}\Big) {\widetilde{\bb}^{(n)}}.
\label{eq:N6Cheb}
\end{equation}

It remains to define the lower block bidiagonal matrix $\tbLn_q$ using
the block Cholesky factorization
\begin{equation}
  \frac{2}{\tau^2} \Big({\bf I}_{zm} - \tilde{\boldsymbol{\mathscr
      P}}_q^{(n)} \Big) = \tbLn_q \tbLnT_q,
  \label{eq:R40}
\end{equation}
calculated with \cite[Algorithm 4.1]{DtB}.  With this matrix we write
\eqref{eq:N6} in the first order algebraic system form
\begin{align}
 \frac{\widetilde \bP_{k+1}^{(n)}-\widetilde \bP_k^{(n)}}{\tau} &= -
 \tbLn_q \thP_k, \quad k = 0, \ldots,
 2n-2, \label{eq:N10}\\ \frac{\thP_k-\thP_{k-1}}{\tau} &=
 \tbLnT_q\widetilde\bP_k^{(n)}, \quad k = 1, \ldots,
 2n-1, \label{eq:N11}
\end{align}
with initial conditions 
\begin{equation}
\widetilde \bP_0^{(n)} =\widetilde{\bb}^{(n)}, \quad \thP_0 +
\thP_{-1} = {\bf 0}.
\label{eq:N12}
\end{equation}
The regularized DtB transformation is
\begin{align}
\widetilde{\bD}_{k}^{\rm Born} & = \widetilde \bD_{0,k} + 
\widetilde \bb^{(n)T} \frac{d}{d \ep} \cT_k \Big(
{\bf I}_{zm} - \frac{\tau^2}{2}
\big[ \tbLn_0 + \ep \big( \tbLn_q - \tbLn_0 \big) \big] 
\big[ \tbLn_0 + \ep \big( \tbLn_q - \tbLn_0\big) \big]^T
\Big) \Big|_{\ep = 0} 
\widetilde \bb^{(n)} \nonumber \\ 
& = \widetilde \bD_{0,k} + \widetilde \bb^{(n)T} 
\delta \widetilde \bP_k^{(n)},
\label{eq:N14}
\end{align}
for $k = 0, 1, \ldots, 2n-1$, where $\widetilde \bD_{0,k}$ and
$\tbLn_0$ are calculated the same way as above, in the reference
medium, {with one important distinction:
The projection matrix $\bZ$ is still computed
from the data corresponding to the medium with (unknown)
reflectivity $q$ and not to the reference medium $q = 0$.
This ensures the consistency of the term 
$\ep \big( \tbLn_q - \tbLn_0 \big)$ in \eqref{eq:N14}.}

{The derivative term in \eqref{eq:N14} is determined by 
the perturbed snapshots $\delta \widetilde \bP_k^{(n)}$, the solutions of}
\begin{align}
\frac{\delta \widetilde \bP_{k+1}^{(n)} - 
\delta \widetilde \bP_k^{(n)}}{\tau} + 
\tbLn_0 \delta \thP_k 
& = - \big( \tbLn_q - \tbLn_0 \big) \thP_{0,k}, 
\quad k = 0, \ldots, 2n-2, 
\label{eq:N15} \\ 
\frac{\delta \thP_{k}-\delta \thP_{k-1}}{\tau} 
- \tbLnT_0 \delta \widetilde \bP_k^{(n)} 
& = \big( \tbLnT_q - \tbLnT_0 \big) 
\widetilde{\bP}_{0,k}^{(n)} , 
\quad k = 1, \ldots, 2n-1, \label{eq:N16}
\end{align}
with homogeneous initial conditions 
\begin{equation}
\delta \widetilde \bP_0^{(n)} = {\bf 0}, \quad 
\delta \thP_0 + \delta \thP_{-1} = {\bf 0}.
\label{eq:N17}
\end{equation}

%=================================================
\subsubsection{Summary of the robust DtB algorithm}
\label{sect:SumAlg}
We now summarize all the steps in the following algorithm that approximates 
the {regularized} DtB transformation \eqref{eq:DtBmap}:

\vspace{0.05in}
\begin{alg}[Robust DtB algorithm] 
\label{alg:DtB}

\textbf{Input:} Data $(\bD_k)_{0 \le k \le 2n-1}.$ 

\textbf{Processing steps:}

\vspace{0.05in}
\begin{enumerate}
\itemsep 0.05in
\item Calculate $\bM$ and $\bS$ using \eqref{eq:14} and \eqref{eq:14p}.
\item Calculate the eigenvalue decomposition \eqref{eq:eigM} and
  define the projection matrix $\bZ$ as in \eqref{eq:defZ}.
\item Generate the data $(\bD_{0,k})_{0 \le k \le 2n-1}$ using
  \eqref{eq:1}--\eqref{eq:3} in the reference medium with reflectivity
  zero. Calculate $\bM_0$ and $\bS_0$ from these data using
  \eqref{eq:14} and \eqref{eq:14p}.
\item Calculate the matrices $ \bSig^2 = \bZ^T \bM \bZ$ and $\cbS =
  \bZ^T \bS \bZ. $ Calculate also their analogues in the reference
  medium $ \bSig_0^2 = \bZ^T \bM_0 \bZ$ and $\cbS_0 = \bZ^T \bS_0
  \bZ.$ Note that $\bSig$ is diagonal but $\bSig_0$ is
  not. Nevertheless, $\bSig_0$ is symmetric and positive definite.
\item Calculate the reduced order model propagator
  $\tilde{\boldsymbol{\mathscr P}}_q^{(n)}$ as the block tridiagonal
  matrix returned by the block Lanczos algorithm for the matrix
  $\bSig^{-1} \cbS \bSig^{-1}$  {and the initial block}
  $\bSig  \widetilde \bE_1 = \bSig \bZ^T \bE_1$. This algorithm produces 
  the orthogonal matrix $\boldsymbol{\mathcal U}$ which defines
  $\widetilde \bb^{(n)}$ as in \eqref{eq:N6p} and
  \eqref{eq:N6ppp}. 
  Similarly, calculate $\tilde{\boldsymbol{\mathscr P}}_0^{(n)}$ 
  using the block Lanczos algorithm for the matrix
  $\bSig^{-1}_0 \cbS_0 \bSig^{-1}_0$.
\item Calculate $\widetilde \bD_{0,k} = \widetilde \bb^{(n)T} \cT_k
  \Big(\tilde{\boldsymbol{\mathscr P}}_0^{(n)}\Big) \widetilde \bb^{(n)}$,
  for $k = 0, \ldots, 2n-1.$
\item Calculate $\Big(\widetilde{\bD}_k^{\rm Born}\Big)_{0 \le k \le  2n-1}$ 
as in \eqref{eq:N14}, with $\delta \widetilde\bP_k^{(n)}$
obtained from the time stepping scheme \eqref{eq:N15}--\eqref{eq:N17}.
\end{enumerate}

\textbf{Output:} The transformed data matrices
$\Big(\widetilde{\bD}_k^{\rm Born}\Big)_{0 \le k \le 2n-1}$.
\end{alg}

\vspace{0.05in} 
The linear algebraic nature of  the operations  in 
Algorithm \ref{alg:DtB} makes it versatile and  useful as a black-box 
active array data processing tool for sound, electromagnetic and elastic waves, 
as we explain in the next three sections.

%=================================================
\section{Application to sonar arrays}
\label{sect:Son}

In this section we show that the acoustic wave equation can be put in
the form \eqref{eq:1}--\eqref{eq:2} and that data gathered by active
sonar arrays can be modeled by \eqref{eq:3}.
\subsection{The sound wave}
\label{sect:Son.1}
Consider the wave equation for the acoustic pressure $p(t,\bx)$ in a
stationary and isotropic medium with wave speed $c(\bx)$ and acoustic
impedance $\sigma(\bx)$,
\begin{equation}
\partial_t^2 p(t,\bx) + A p(t,\bx) = \partial_t f(t) \delta(\bx -
\bx_s), 
\label{eq:S1}
\qquad A p(t,\bx)= - \sigma(\bx) c(\bx) \nabla \cdot \Big[
  \frac{c(\bx)}{\sigma(\bx)} \nabla p(t,\bx) \Big],
\end{equation}
for $t \in \RR$ and $\bx \in \Omega$, where we note that $\sigma(\bx)
c(\bx)$ is the bulk modulus and $\sigma(\bx)/c(\bx)$ is the mass
density. The wave is generated by the source at $\bx_s$, which emits
the pulse $f(t)$ that is compactly supported around $t = 0$. Prior to
the excitation the medium is at equilibrium
\begin{equation}
p(t,\bx) = 0, \quad t \ll 0.
\label{eq:S1p}
\end{equation}
As explained in the previous section, since the wave equation is causal and the wave speed is finite,
we can restrict $\bx$ to a compact cube $\Omega_c \subset \Omega$ of
side length larger than $2 n \tau \max{c(\bx)}$ and with boundary
$\partial \Omega_c$ given by the union of two sets: $\partial
\Omega_c^{\rm ac}$ contained in $\partial \Omega$, called the
accessible boundary, and the inaccessible boundary $\partial
\Omega_c^{\rm inac}$ contained in $\Omega$. As in \cite{DtB} we model the accessible
boundary as sound hard,
\begin{equation}
\partial_{x_d} p(t,\bx) = 0, \quad \bx \in \partial \Omega_c^{\rm ac},
\label{eq:S2}
\end{equation}
but other homogeneous boundary conditions can be considered, as well.
The center of $\partial \Omega_c^{\rm ac}$ is assumed just below the
center of the array, which has aperture $a \ll 2 n \tau \max{c(\bx)}$.
At the inaccessible boundary we set
\begin{equation}
p(t,\bx) = 0, \quad \bx \in \partial \Omega_c^{\rm inac},
\label{eq:S3}
\end{equation}
without affecting the wave over the duration $|t| < 2 n \tau$ of the
measurements.

%=================================================
\subsection{Data model}
\label{sect:DMson}
The operator $A$ defined in \eqref{eq:S1} for $\bx \in \Omega_c$, with
the homogeneous boundary conditions \eqref{eq:S2}--\eqref{eq:S3}, is
positive definite and self-adjoint in the Hilbert space $ L^2\Big(
\Omega_c,\frac{1}{\sigma(\bx) c(\bx)} {\rm d} \bx \Big)$ with weighted
inner product
\begin{equation}
\left<\varphi,\psi\right>_{\frac{1}{\sigma c}} = \int_{\Omega_c} {\rm
  d} \bx \, \frac{\varphi(\bx) \psi(\bx)}{\sigma(\bx) c(\bx)}.
\label{eq:S4}
\end{equation}
We define the data for the $s^{\rm th}$ source excitation using the
even extension in time of the pressure wave measured at the receiver
locations $\bx_r$, for $r = 1, \ldots, m_a$,
\begin{equation}
D_k^{(r,s)} = p^{e}(t_k,\bx_r) = p(t_k,\bx_r) + p(-t_k,\bx_r), \quad
t_k = k \tau, ~ ~ k = 0, \ldots, 2n-1.
\label{eq:S4p}
\end{equation}
This even extension gives the homogeneous initial condition  $\partial_t
p^{(e)}(0,\bx) = 0$, as in \eqref{eq:2}.

The wave $p^{e}(t_k,\bx_r)$ can be written more explicitly, if we make
the convenient assumption that the pulse has a real valued and
non-negative Fourier transform $\hat f(\omega)$. Then, as shown in
\cite[Section 2.1]{DtB},
\begin{equation}
p^{e}(t,\bx) = \cos(t \sqrt{A}) \hat f(\sqrt{A}) \delta(\bx-\bx_s).
\label{eq:S5}
\end{equation}
Moreover, if we let 
\begin{equation}
b^{(s)}(\bx) = \sqrt{\frac{\sigma(\bx_s) c(\bx_s)}{\sigma(\bx)
    c(\bx)}} \Big[\hat f(\sqrt{A})\Big]^{1/2} \delta(\bx-\bx_s),
\label{eq:S6}
\end{equation}
we can rewrite \eqref{eq:S4p} in the symmetric form 
\begin{equation}
D_k^{(r,s)} = \left< \sqrt{\sigma c}\,  b^{(r)}, \cos( t_k \sqrt{A})
\sqrt{\sigma c} \, b^{(s)}\right>_{\frac{1}{\sigma c}}, \qquad k = 0, \ldots, 2n-1.
\label{eq:S7}
\end{equation}

\vspace{0.05in}
\begin{rem}
\label{rem:sensors}
The function $b^{(s)}(\bx)$ is localized near the source location
$\bx_s$, so we think of it as a sensor indicator function. The data
model can be interpreted as the wave $\cos( t \sqrt{A}) \sqrt{\sigma
  c} \, b^{(s)}$ generated by the source modeled by $b^{(s)}$,
measured at time $t_k$ by the receiver modeled by $b^{(r)}$.
\end{rem}

%=================================================
\subsection{The reflectivity model}
The joint estimation of $c(\bx)$ and $\sigma(\bx)$ from the data
\eqref{eq:S7} is difficult, especially for the small array aperture
considered in our setting. Moreover, these coefficients play a
different role in the wave propagation process: While the smooth part
of the velocity determines the kinematics of the wave i.e., the travel
times, the variations of the acoustic impedance determine the dynamics
of the wave i.e., the reflections.

It is shown in \cite[Section 1]{BeylkinBurridge} that if the array is
small, the main contribution to the single scattered wave field, the
Born approximation, is determined by the variations of
\begin{equation}
q(\bx) = \ln \sigma(\bx),
\label{eq:S8}
\end{equation}
{which we refer to as the reflectivity.} We assume henceforth 
that we know $c(\bx)$, although in practice we can only know its smooth part.
Depending on the application, this could be a constant or it could be
a function estimated independently, with a velocity estimation method
like in
\cite{uhlmann2001travel,uhlmann2016inverse,stolk2002smooth,shen2003differential}.
If the medium has constant mass density, then the variations of
$\sigma(\bx)$ determine the variations of $c(\bx)$. Otherwise, the
rough variations of $c(\bx)$ alone, appear in the expression of the
Born approximation multiplied by $\sin^2 \vartheta$, where $\vartheta$
is the angle between the outgoing and incoming rays connecting the
source and receiver to a point of reflection \cite[Section
  1]{BeylkinBurridge}. When the array is small and the variations in
the medium are sufficiently far from the surface, the angle
$\vartheta$ is small, so we expect that the leading contribution
comes from the reflectivity \eqref{eq:S8}.

%=================================================
\subsection{The Liouville transform}
To write the problem in the  form \eqref{eq:1}, let us
introduce the wave
\begin{equation}
p^{(s)}(t,\bx) = \cos (t \sqrt{A}) \sqrt{\sigma(\bx) c(\bx)}\, b^{(s)}(\bx),
\label{eq:S9}
\end{equation}
indexed by the source.  This is the pressure field in the first order
system of acoustic wave equations
\begin{equation}
\partial_t \begin{pmatrix} p^{(s)}(t,\bx)
  \\ \bu^{(s)}(t,\bx) \end{pmatrix} =
\begin{pmatrix}0 & - \sigma(\bx) c(\bx) \nabla \cdot \\
-\frac{c(\bx)}{\sigma(\bx)} \nabla & 0 \end{pmatrix} \begin{pmatrix}
  p^{(s)}(t,\bx) \\ \bu^{(s)}(t,\bx) \end{pmatrix}, \quad \bx \in
\Omega_c, \quad t > 0,
\label{eq:dtb1}
\end{equation}
with boundary conditions 
\begin{align}
{\bf e}_d \cdot \bu^{(s)}(t,\bx) = 0, \quad \bx \in \partial
\Omega_c^{\rm ac} \label{eq:S10} \quad \mbox{and} \quad p^{(s)}(t,\bx) = 0, \quad
\bx \in \partial \Omega_c^{\rm inac},
\end{align}
and initial conditions
\begin{equation}
p^{(s)}(0,\bx) = \sqrt{\sigma(\bx) c(\bx)} b^{(s)}(\bx), \quad
\bu^{(s)}(0,\bx) = {\bf 0}.\label{eq:S12}
\end{equation}
Here $\bu^{(s)}(t,\bx)$ is the dual wave, the displacement velocity,
and ${\bf e}_d$ is the unit vector along the $x_d$ axis in $\RR^{d}$.

Let us also use a Liouville transformation to define the primary wave
\begin{equation}
P^{(s)}(t,\bx) = \frac{ p^{(s)}(t,\bx)}{\sqrt{\sigma(\bx) c(\bx)}},
\label{eq:S13}
\end{equation}
which is scalar valued $(\dP = 1)$ and the dual wave
\begin{equation}
\hat {\bP}^{(s)}(t,\bx) = - \sqrt{\frac{\sigma(\bx)}{c(\bx)}}
\bu^{(s)}(t,\bx),
\label{eq:S14}
\end{equation}
which is vector valued. Substituting in \eqref{eq:dtb1}--\eqref{eq:S12}, we obtain that these
satisfy the first order system
\begin{equation}
\partial_t \begin{pmatrix} P^{(s)}(t,\bx) \\ \widehat
  {\bP}^{(s)}(t,\bx) \end{pmatrix} = \begin{pmatrix} 0 & - L_q
  \\ L_q^T & 0 \end{pmatrix} \begin{pmatrix} P^{(s)}(t,\bx) \\\widehat
  {\bP}^{(s)}(t,\bx) \end{pmatrix}, \quad \bx \in \Omega_c, \quad t >
0,
\label{eq:S15}
\end{equation}
with boundary conditions
\begin{align}
{\bf e}_d \cdot \widehat \bP^{(s)}(t,\bx) = 0, \quad \bx \in \partial
\Omega_c^{\rm ac} \label{eq:S16} \quad \mbox{and} \quad P^{(s)}(t,\bx)
= 0, \quad \bx \in \partial \Omega_c^{\rm inac},
\end{align}
for $t > 0$, and initial conditions
\begin{equation}
P^{(s)}(0,\bx) = b^{(s)}(\bx), \quad \widehat \bP^{(s)}(0,\bx) = {\bf
  0}, \qquad \bx \in \Omega_c.\label{eq:S18}
\end{equation}
The first order operator $L_q$ is given by 
\begin{equation}
L_q\widehat \bP^{(s)}(t,\bx) = - \sqrt{c(\bx)} \nabla \cdot \Big[
  \sqrt{c(\bx)}\widehat \bP^{(s)}(t,\bx) \Big] + \frac{c(\bx)}{2}
\nabla q(\bx) \cdot \widehat \bP^{(s)}(t,\bx),
\label{eq:S19}
\end{equation}
and its adjoint with respect to the Euclidian inner product is 
\begin{equation}
L_q^T P^{(s)}(t,\bx) = \sqrt{c(\bx)} \nabla \Big[ \sqrt{c(\bx)}
  P^{(s)}(t,\bx) \Big] + \frac{c(\bx)}{2} \nabla q(\bx)
P^{(s)}(t,\bx).
\label{eq:S20}
\end{equation}
These are affine with respect to the reflectivity $q(\bx)$.

The generic model \eqref{eq:1}--\eqref{eq:2} follows once we write
\eqref{eq:S15}--\eqref{eq:S18} as a second order wave equation for
$P^{(s)}(t,\bx)$. Moreover, the data \eqref{eq:S7} can be rewritten as
\begin{equation}
D_k^{(r,s)} = \int_{\Omega_c} {\rm d} \bx \, b^{(r)}(\bx)
P^{(s)}(t_k,\bx), \qquad k = 0, \ldots, 2n-1,
\label{eq:S21}
\end{equation}
where the integral over $\Omega_c$ can be extended to the whole domain
$\Omega$, because the wave vanishes in $\Omega \setminus \Omega_c$
over the duration of the measurements. This is the same as
\eqref{eq:3}.

%=================================================
\section{Application to electromagnetic waves}
\label{sect:Max}
In this section we show that Maxwell's equations can be put in the form
\eqref{eq:1}--\eqref{eq:2} and that the array data can be modeled by
\eqref{eq:3}.

%=================================================
\subsection{The electric field}
Consider the wave equation for the electric field
${\boldsymbol{\mathcal E}}^{(s)}(t,\bx)$ in an isotropic, lossless and
time independent medium with electric permittivity $\vep(\bx)$ and
magnetic permeability $\mu(\bx)$,
\begin{equation}
\partial_t^2 {\boldsymbol{\mathcal E}}^{(s)}(t,\bx) +
\frac{1}{\vep(\bx)} \nabla \times \Big[\frac{1}{\mu(\bx)} \nabla
  \times {\boldsymbol{\mathcal E}}^{(s)}(t,\bx)\Big] = 0, \qquad t >
0, \quad \bx \in \Omega \subset \RR^3,
\label{eq:E1}
\end{equation}
with initial conditions 
\begin{equation}
{\boldsymbol{\mathcal E}}^{(s)}(0,\bx) =
\frac{\bb^{(s)}(\bx)}{\sqrt{\vep(\bx)}}, \quad
\partial_t{\boldsymbol{\mathcal E}}^{(s)}(0,\bx) = {\bf 0}, \qquad \bx
\in \Omega.
\label{eq:E2}
\end{equation}
These can be derived as in the sonar case, starting from the equation
with a forcing term due to a source in the array, and then taking the
even extension in time. The initial condition is defined by the
``sensor function'' $\bb^{(s)}$, which is localized near the source
and is defined with a similar procedure to that in section
\ref{sect:DMson}. We suppose that it satisfies
\begin{equation}
\nabla \cdot \Big[ \vep(\bx) {\boldsymbol{\mathcal E}}^{(s)}(0,\bx) \Big]  = 0,
\label{eq:E3}
\end{equation}
so that the electric displacement $\vep(\bx) {\boldsymbol{\mathcal
    E}}(t,\bx)$ remains divergence free at all times.

As in the previous sections, we restrict the domain $\Omega$ to a
compact cube $\Omega_c$ with boundary $\partial \Omega_c = \partial
\Omega_c^{\rm ac} \cup \partial \Omega_c^{\rm inac}$ consisting of
the accessible boundary $\partial\Omega_c^{\rm ac} \subset \partial
\Omega$ and the inaccessible boundary $\partial\Omega_c^{\rm inac}$
contained in $\Omega$. For simplicity, we model both boundaries as
perfectly conducting,
\begin{equation}
\boldsymbol{\nu}(\bx) \times {\boldsymbol{\mathcal E}}^{(s)}(t,\bx) =
           {\bf 0}, \qquad \bx \in \partial \Omega_c, ~~ t > 0,
\label{eq:E4}
\end{equation}
where $\boldsymbol{\nu}(\bx)$ is the unit outer normal at $\partial
\Omega_c$. Other homogeneous boundary conditions can be considered, as
well.

%=================================================
\subsection{The reflectivity}
The electric permittivity and magnetic permeability define the wave
speed $c(\bx)$ and the impedance $\sigma(\bx)$,
\begin{equation}
c(\bx) = \frac{1}{\sqrt{\mu(\bx) \vep(\bx)}}, \qquad \sigma(\bx) =
\sqrt{\frac{\mu(\bx)}{\vep(\bx)}}.
\label{eq:E5}
\end{equation}
Like in sonar, we assume that the wave speed is known, and define the
unknown reflectivity by
\begin{equation}
q(\bx) = \ln \sigma(\bx).
\label{eq:E6}
\end{equation}
Again, in practice, only the smooth part of the wave speed is known.
The rough variations of $c(\bx)$ are captured by $q(\bx)$ when the
magnetic permeability is constant. Even when $\mu(\bx)$ varies, if the
array is small, it is not $c(\bx)$ alone that plays the leading role
in the reflection, but $\sigma(\bx) = c(\bx) \mu(\bx)$ or,
equivalently, the reflectivity $q(\bx)$.

%=================================================
\subsection{The data model}
We can rewrite the initial condition \eqref{eq:E2} in terms of
$\sigma(\bx)$ and $c(\bx)$, using \eqref{eq:E5},
\begin{equation}
{\boldsymbol{\mathcal E}}^{(s)}(0,\bx) = \sqrt{\sigma(\bx) c(\bx)} \bb^{(s)}(\bx),
\label{eq:E7}
\end{equation}
and proceed as in section \ref{sect:DMson} (see Remark
\ref{rem:sensors}) to write the data model in the symmetric form
\begin{equation}
D_k^{(r,s)} = \left< \sqrt{\sigma c} \, \bb^{(r)},
{\boldsymbol{\mathcal E}}^{(s)}(t_k,\bx) \right>_{\frac{1}{\sigma c}}
= \left< \sqrt{\sigma c} \, \bb^{(r)}, \cos(t_k \sqrt{A}) \sqrt{\sigma
  c} \, \bb^{(s)} \right>_{\frac{1}{\sigma c}},
\label{eq:E8}
\end{equation}
for $t_k = k \tau$, $k = 0, \ldots, 2n-1$ and $r,s = 1, \ldots, m$.
Here $A$ is the operator
\begin{equation}
A {\boldsymbol{\mathcal E}}^{(s)}(t,\bx) = \sigma(\bx) c(\bx) \nabla
\times \Big[\frac{c(\bx)}{\sigma(\bx)} \nabla \times
  {\boldsymbol{\mathcal E}}^{(s)}(t,\bx)\Big],
\label{eq:E9}
\end{equation}
defined on functions satisfying the boundary conditions \eqref{eq:E4},
which is positive definite and self-adjoint with respect to the
weighted inner product
\begin{equation}
\left< \boldsymbol{\varphi}, \boldsymbol{\psi}
\right>_{{\frac{1}{\sigma c}}} = \int_{\Omega_c} {\rm d} \bx \,
\frac{\boldsymbol{\varphi}^T(\bx) \boldsymbol{\psi}(\bx)}{\sigma(\bx)
  c(\bx)}.
\end{equation}

Note that the ``sensor functions'' $\bb^{(s)}(\bx)$ are now vector
valued, and the index $s$ stands not only for the sensor location but
also for the polarization of the wave. It is sufficient to prescribe
and measure the tangential components of ${\boldsymbol{\mathcal
    E}}^{(s)}(t,\bx)$, so we have two functions per sensor location,
meaning that $m = 2 m_a$.

%=================================================
\subsection{The Liouville transformation}
The first order system formulation of \eqref{eq:E1}--\eqref{eq:E2} is
\begin{equation}
\partial_t \begin{pmatrix} 
{\boldsymbol{\mathcal E}}^{(s)}(t, \bx) \\ 
{\boldsymbol{\mathcal H}}^{(s)}(t, \bx)
\end{pmatrix} = 
\begin{pmatrix} 0 & \sigma(\bx) c(\bx) \nabla \times \\
- \frac{c(\bx)}{\sigma(\bx)} \nabla \times & 0 
\end{pmatrix}
\begin{pmatrix} 
{\boldsymbol{\mathcal E}}^{(s)}(t, \bx) \\ 
{\boldsymbol{\mathcal H}}^{(s)}(t, \bx)
\end{pmatrix}, 
\qquad t > 0, ~ ~ \bx \in \Omega_c,
\label{eq:E10} 
\end{equation}
with initial conditions
\begin{equation}
{\boldsymbol{\mathcal E}}^{(s)}(0,\bx) = \sqrt{\sigma(\bx) c(\bx)}
\bb^{(s)}(\bx), \qquad {\boldsymbol{\mathcal H}}^{(s)}(0,\bx) = {\bf
  0}, \qquad \bx \in \Omega_c,
\label{eq:E11}
\end{equation}
where ${\boldsymbol{\mathcal H}}^{(s)}(t,\bx)$ is the magnetic field.
We now use a Liouville transformation of the electric and magnetic
fields to derive the generic wave equation \eqref{eq:1}.

The primary wave is the vector valued field ($\dP = 3$)
\begin{equation}
\bP^{(s)}(t,\bx) = \frac{{\boldsymbol{\mathcal
      E}}^{(s)}(t,\bx)}{\sqrt{\sigma(\bx)c(\bx)}},
\label{eq:E12}
\end{equation}
and the dual wave is 
\begin{equation}
\widehat \bP^{(s)}(t,\bx) =
\sqrt{\frac{\sigma(\bx)}{c(\bx)}}{\boldsymbol{\mathcal H}}^{(s)}(t,\bx).
\label{eq:E13}
\end{equation}
Substituting these definitions in \eqref{eq:E4},
\eqref{eq:E10}--\eqref{eq:E11} we obtain that
\begin{equation}
\partial_t \begin{pmatrix} \bP^{(s)}(t,\bx) \\ \widehat
  {\bP}^{(s)}(t,\bx) \end{pmatrix} = \begin{pmatrix} 0 & - L_q
  \\ L_q^T & 0 \end{pmatrix} \begin{pmatrix} \bP^{(s)}(t,\bx)
  \\\widehat {\bP}^{(s)}(t,\bx) \end{pmatrix}, \quad \bx \in \Omega_c,
\quad t > 0,
\label{eq:E14}
\end{equation}
with boundary conditions 
\begin{align}
\boldsymbol{\nu}(\bx) \times \bP^{(s)}(t,\bx) = 0, \quad \bx \in
\partial \Omega_c, \quad t > 0,
\label{eq:E15}
\end{align}
and initial conditions
\begin{equation}
\bP^{(s)}(0,\bx) = \bb^{(s)}(\bx), \qquad \widehat \bP^{(s)}(0,\bx) =
   {\bf 0}, \quad \bx \in \Omega_c.
\label{eq:E16}
\end{equation}
The first order operator $L_q$ is given by 
\begin{equation}
L_q \widehat \bP^{(s)}(t,\bx) = -\sqrt{c(\bx)}\nabla \times \Big[
  \sqrt{c(\bx)} \widehat \bP^{(s)}(t,\bx)\Big] +
\frac{c(\bx)}{2} \nabla q(\bx) \times \widehat
\bP^{(s)}(t,\bx),
\label{eq:E17}
\end{equation}
and $L_q^T$ is its adjoint with respect to the Euclidian inner product,
\begin{equation}
L_q^T \bP^{(s)}(t,\bx) = -\sqrt{c(\bx)}\nabla \times \Big[
  \sqrt{c(\bx)} \bP^{(s)}(t,\bx)\Big] - \frac{c(\bx)}{2} \nabla q(\bx)
\times \bP^{(s)}(t,\bx).
\label{eq:E18}
\end{equation}

The generic model \eqref{eq:1}--\eqref{eq:2} follows by writing
\eqref{eq:E17}--\eqref{eq:E18} as a second order wave equation for the
primary wave $\bP^{(s)}(t,\bx)$. Moreover, we obtain from \eqref{eq:E8} 
and \eqref{eq:E12} the data model 
\begin{equation}
D_k^{(r,s)} = \int_{\Omega_c} {\rm d} \bx \, \big(\bb^{(r)}(\bx)\big)^T 
\bP^{(s)}(t_k,\bx), \qquad k = 0, \ldots, 2n-1,
\label{eq:E19}
\end{equation}
which is the same as \eqref{eq:3}, because the wave is zero in $\Omega
\setminus \Omega_c$.

%=================================================
\section{Application to elastic waves}
\label{sect:EL}
In this section we consider elastic waves in an isotropic, time
independent medium with mass density $\rho(\bx)$ and Lam\'{e}
parameters $\lambda(\bx)$ and $\mu(\bx)$.  For simplicity, we derive
the wave model \eqref{eq:1}--\eqref{eq:3} in two dimensions ($d =
2$). As in the previous sections, we restrict the wave to the compact
cube $\Omega_c \subset \Omega$, with large enough side length so that
the waves are not affected over the duration $2 n \tau$ by the
boundary conditions at the inaccessible part $\partial \Omega_c^{\rm
  inac}$ of $ \partial \Omega_c$.

%=================================================
\subsection{The elastic wave equation and the data model}
We begin with Newton's second law
\begin{equation}
\rho(\bx) \partial_t^2 u_i^{(s)}(t,\bx) = \sum_{j=1}^2 \partial_{x_j}
\St_{ij}^{(s)}(t,\bx), \quad i = 1,2, \quad t > 0, \quad \bx \in
\Omega_c,
\label{eq:EL2}
\end{equation}
where $\bu^{(s)}(t,\bx) = \big(u_j^{(s)}(t,\bx)\big)_{1 \le j \le 2}$
is the displacement vector and $\Big(\St_{i,j}^{(s)}\Big)_{1 \le i, j \le 2}$
is the {symmetric} stress tensor, with components
\begin{equation}
\St_{ij}^{(s)}(t,\bx) = \lambda(\bx) \nabla \cdot \bu^{(s)}(t,\bx)
\delta_{ij} + \mu(\bx) 
{\left[ \partial_{x_j} u_i^{(s)}(t,\bx) + 
\partial_{x_i} u_j^{(s)}(t,\bx)\right], }
\quad i,j = 1,2.
\label{eq:EL1}
\end{equation}
The index $(s)$ stands for the $s^{\rm th}$ source excitation, which
is defined below via initial conditions.

Now let us introduce the velocity
\begin{equation}
\bv^{(s)}(t,\bx) = \partial_t \bu^{(s)}(t,\bx),
\label{eq:EL3}
\end{equation}
with components $v_j^{(s)}(t,\bx)$, for $j = 1, 2$, and organize the five
unknowns $v_1^{(s)},v_2^{(s)}, \St^{(s)}_{11}, \St^{(s)}_{22}, \St^{(s)}_{12}$ in the vector
$\begin{pmatrix} \bv^{(s)} \\ \bSt^{(s)} \end{pmatrix} $. We obtain from
\eqref{eq:EL2} the first order system
\begin{align}
\partial_t \begin{pmatrix} \bv^{(s)}(t,\bx) \\ \bSt^{(s)}(t,\bx) \end{pmatrix}
= \begin{pmatrix} {\bf 0} & - \frac{1}{\rho} \bcD^T \\ \bK(\bx) \bcD &
  {\bf 0} \end{pmatrix}
\begin{pmatrix} \bv^{(s)}(t,\bx) \\ \bSt^{(s)}(t,\bx) \end{pmatrix}, 
\quad t > 0, \quad \bx \in \Omega_c,
\label{eq:EL4}
\end{align}
with coefficient matrix
\begin{equation}
\bK(\bx) = \begin{pmatrix} \lambda(\bx) + 2 \mu(\bx) & \lambda(\bx) &
  0 \\ \lambda(\bx) & \lambda(\bx) + 2 \mu(\bx) & 0 \\ 0 & 0 &
  \mu(\bx)
\end{pmatrix} 
\label{eq:2.n0}
\end{equation}
and differential operators  
\begin{equation}
-\bcD^T = \begin{pmatrix} \partial_{x_1} & 0 & \partial_{x_2} \\ 0 &
  \partial_{x_2} & \partial_{x_1} \end{pmatrix}, \qquad \bcD
= \begin{pmatrix} \partial_{x_1} & 0 \\0 & \partial_{x_2}
  \\ \partial_{x_2} & \partial_{x_1} \end{pmatrix}.
\label{eq:2.n2}
\end{equation}
This system is endowed with traction free (Neumann) boundary conditions at the accessible boundary 
with outer normal   $\boldsymbol{\nu} = (\nu_1,\nu_2) $,
\begin{equation}
\begin{pmatrix} 
\nu_1(\bx) & 0 & \nu_2(\bx)\\
0 & \nu_2(\bx) & \nu_1(\bx) 
\end{pmatrix} \bSt^{(s)}(t,\bx) = 0,  \quad  \bx \in \partial \Omega_c^{\rm ac}, \quad  t > 0, 
\label{eq:EL5N}
\end{equation}
the homogeneous Dirichlet conditions at the inaccessible boundary
\begin{equation}
\bv^{(s)} = {\bf 0}, \quad  \bx \in \partial \Omega_c^{\rm inac},  \quad t > 0,
\label{eq:EL5}
\end{equation}
and the initial conditions
\begin{equation}
\bv^{(s)}(0,\bx) = \sqrt{\rho(\bx)} \bb^{(s)}(\bx),
\quad \bSt^{(s)}(0,\bx) = {\bf 0}, \quad \bx \in \Omega_c.
\label{eq:EL6}
\end{equation}
These initial conditions can be defined as in section \ref{sect:Son}, by
considering {an external force term in \ref{eq:EL2} due to an array sensor and assuming} the even extension in time of the velocity, which is
consistent with setting $\partial_t \bv^{(s)}(0,\bx) = {\bf 0}$ or,
equivalently, $\bSt^{(s)}(0,\bx) = {\bf 0}.$ { We note that there are two independent stress-strain state excitations at each sensor location.} Moreover, the initial
velocity is localized near the source location, and it is proportional
to the ``sensor function'' $\bb^{(s)}(\bx)$, defined with a similar
procedure as in section \ref{sect:DMson}. Recall also Remark
\ref{rem:sensors}.

The data are defined by
\begin{equation} 
D_{k}^{(r,s)} = \int_{\Omega_c} {\rm d} \bx \, 
\big(\bb^{(r)}(\bx) \big)^T\sqrt{\rho(\bx) }\,\bv^{(s)}(t_k,\bx),
\label{eq:ELDat}
\end{equation}
for $t_k = k \tau$, $k = 0,\ldots, 2n-1$ and $r,s = 1, \ldots,
m$. Since $\bb^{(s)}(\bx) \in \RR^2$, we can have two different
excitations per source location and two components of the velocity
measured at the receivers, meaning that $m = 2 m_a$.

%=================================================
\subsection{The wave speeds and impedances}
\label{sect:EL2}
The density and the Lam\'{e} parameters define the wave speed and
impedance for the pressure wave
\begin{equation}
c_p(\bx) = \sqrt{\frac{\lambda(\bx)+2 \mu(\bx)}{\rho(\bx)}}, \quad
\sigma_p(\bx) = c_p(\bx) \rho(\bx),
\label{eq:2.17}
\end{equation}
and the wave speed and impedance for the shear wave
\begin{equation}
c_s(\bx)= \sqrt{\frac{\mu(\bx)}{\rho(\bx)}}, \quad \sigma_s(\bx) =
c_s(\bx) \rho(\bx).
\label{eq:2.18}
\end{equation}
These definitions give the following identities
\begin{align*}
\frac{1}{\rho(\bx)} &= \frac{c_p(\bx)}{\sigma_p(\bx)} =
\frac{c_s(\bx)}{\sigma_s(\bx)} ,\\ \lambda(\bx) + 2 \mu(\bx) & =
c_p(\bx) \sigma_p(\bx), \\ \mu(\bx) &= c_s(\bx) \sigma_s(\bx),
\\ \sigma_s(\bx) &= \sigma_p(\bx) \frac{c_s(\bx)}{c_p(\bx)}. 
\end{align*}

Let us introduce the dimensionless parameter
\begin{equation}
\gamma(\bx) = \left(\frac{c_s(\bx)}{c_p(\bx)}\right)^2 =
\frac{\mu(\bx)}{\lambda(\bx) + 2 \mu(\bx)} < 1.
\label{eq:2.19b}
\end{equation}
Then, we have 
\begin{align}
\lambda(\bx) &= c_p(\bx) \sigma_p(\bx) - 2 \mu(\bx) = c_p(\bx)
\sigma_p(\bx) (1-2 \gamma(\bx)) \label{eq:2.19c}, \\\mu(\bx) &=
\sigma_p(\bx) \frac{c_s^2(\bx)}{c_p(\bx)} = c_p(\bx) \sigma_p(\bx)
\gamma(\bx),
\end{align}
and the coefficient matrix \eqref{eq:2.n0} can be rewritten as 
\begin{equation} 
\bK(\bx) = c_p(\bx) \sigma_p(\bx) \bGamma(\bx), \quad \bGamma(\bx)
= \begin{pmatrix} 1 & 1 - 2 \gamma(\bx) & 0 \\ 1-2 \gamma(\bx) & 1 & 0
  \\ 0 & 0 & \gamma(\bx) \end{pmatrix}.
\label{eq:2N3}
\end{equation}
The first order system \eqref{eq:EL4} becomes
\begin{align}
\partial_t \begin{pmatrix} \bv^{(s)}(t,\bx)
  \\ \bSt^{(s)}(t,\bx) \end{pmatrix} = \begin{pmatrix}
  \frac{c_p(\bx)}{\sigma_p(\bx)} {\bf I}_2 & {\bf 0} \\ {\bf 0} &
  c_p(\bx) \sigma_p(\bx) \bGamma(\bx) \end{pmatrix} \begin{pmatrix}
  {\bf 0} & - \bcD^T \\ \bcD & {\bf 0} \end{pmatrix}
\begin{pmatrix} \bv^{(s)}(\bx) \\ \bSt^{(s)}(\bx) \end{pmatrix}, \quad t > 0, \quad 
\bx \in \Omega_c,
\label{eq:2NS1}
\end{align}
where ${\bf I}_2$ is the $2\times 2$ identity. At $\partial \Omega_c$
we have the boundary conditions \eqref{eq:EL5N}--\eqref{eq:EL5} and we rewrite the
initial state in \eqref{eq:EL6} as
\begin{equation}
\bv_o^{(s)}(\bx) = \sqrt{\sigma_p(\bx) c_p(\bx)} \bb^{(s)}(\bx).
\label{eq:EL10}
\end{equation}
Moreover, the data \eqref{eq:ELDat} become, for $t_k = k \tau$, $k = 0,\ldots, 2n-1$ and $r,s = 1, \ldots,
m$,
\begin{equation} 
D_{k}^{(r,s)} = \int_{\Omega_c} {\rm d} \bx \, \big(\bb^{(r)}(\bx)\big)^T
\sqrt{\frac{\sigma_p(\bx)}{c_p(\bx)}}
\bv^{(s)}(t_k,\bx). \label{eq:ELDat1}
\end{equation}
%=================================================
\subsection{The Liouville transformation}

Note that the symmetric matrix $\bGamma(\bx)$ defined in
\eqref{eq:2N3} has the eigenvalue decomposition
\begin{equation}
\bGamma(\bx) = \begin{pmatrix} \frac{1}{\sqrt{2}} &  \frac{1}{\sqrt{2}} & 0 \\
\frac{1}{\sqrt{2}} & - \frac{1}{\sqrt{2}} & 0 \\
0 & 0 & 1 \end{pmatrix} 
\begin{pmatrix} 2(1-\gamma(\bx)) & 0 & 0 \\ 0 & 2 \gamma(\bx) & 0 
\\ 0 & 0 & \gamma(\bx) \end{pmatrix}
\begin{pmatrix} \frac{1}{\sqrt{2}} &  \frac{1}{\sqrt{2}} & 0 \\
\frac{1}{\sqrt{2}} & - \frac{1}{\sqrt{2}} & 0 \\ 0 & 0 &
1 \end{pmatrix}\label{eq:eig}, \end{equation} and  it is positive
definite because $\gamma(\bx) < 1$. We use it in the
Liouville transformation that defines the primary wave
\begin{equation}
\bP^{(s)}(t,\bx) = \sqrt{\frac{\sigma_p(\bx)}{c_p(\bx)}}
\bv^{(s)}(t,\vx),
\label{eq:EL11}
\end{equation}
which is vector valued in $\RR^2$ i.e., $\dP = 2$, and the dual wave
\begin{equation}
\widehat \bP^{(s)}(t,\bx)= \frac{1}{\sqrt{\sigma_p(\bx)c_p(\bx)}
}\bGamma^{-1/2}(t,\bx) \bSt^{(s)}(t,\bx).
\label{eq:2N4}
\end{equation}
These waves satisfy the first order system
\begin{align}
\partial_t \begin{pmatrix} \bP^{(s)}(t,\bx)
  \\ \widehat \bP^{(s)}(t,\bx) \end{pmatrix} =
 \begin{pmatrix} 0 & -  L_{q}
\\ L_{q}^T & 0 \end{pmatrix}
\begin{pmatrix} \bP^{(s)}(t,\bx) \\ \widehat \bP^{(s)}(t,\bx) \end{pmatrix}, \quad 
t > 0, \quad \bx \in \Omega,
\label{eq:2N5}
\end{align}
with operators $L_q$ and $L_q^T$ defined by
\begin{align}
L_{q}\widehat \bP^{(s)}(t,\bx) &= \sqrt{c_p(\bx)} \bcD^T
\left(\sqrt{c_p(\bx)}\bGamma^{1/2}(\bx) \widehat \bP^{(s)}(t,\bx)
\right) - \bQQ_{q}(\bx)\bGamma^{1/2}(\bx) \widehat \bP^{(s)}(t,\bx),
\\
L_{q}^T\bP^{(s)}(t,\bx) &= \sqrt{c_p(\bx)} \bGamma^{1/2}(\bx) \bcD
\Big(\sqrt{c_p(\bx)}\bP^{(s)}(t,\bx)\Big) - \bGamma^{1/2}(\bx)
\bQQ_{q}^T(\bx)\bP^{(s)}(t,\bx).
\end{align}
Here we introduced the $2 \times 3$ matrix valued potential
\begin{equation}
\bQQ_{q}(\bx) = \frac{1}{2} c_p(\bx) \left(-\bcD^T q(\bx) {\bf
  I}_3\right) = \frac{c_p(\bx)}{2} \begin{pmatrix} \partial_{x_1}
  q(\bx) & 0 & \partial_{x_2} q(\bx) \\ 0 & \partial_{x_2}
  q(\bx) & \partial_{x_1}q(\bx) \end{pmatrix},
\end{equation}
where ${\bf I}_3$ is the $3 \times 3$ identity. This potential depends
linearly on the reflectivity
\begin{equation}
q(\bx) = \ln \sigma_p(\bx). 
\label{eq:EL15}
\end{equation}
By definition \eqref{eq:2.n2} the operator $\bcD$ is the adjoint of
$\bcD^T$ and the adjoint potential is the $3 \times 2$ matrix
\begin{equation}
\bQQ_{q}^T(\bx) = \frac{c_p(\bx)}{2} \Big(\bcD q(\bx) {\bf I}_2\Big) =
\frac{c_p(\bx)}{2} \begin{pmatrix} \partial_{x_1} q(\bx) & 0 \\0 &
  \partial_{x_2} q(\bx) \\ \partial_{x_2} q(\bx) &
  \partial_{x_1}q(\bx) \end{pmatrix}.
\end{equation}
Thus, the operator $L_q^T$ is the adjoint of $L_q$.

The model \eqref{eq:1}--\eqref{eq:2} follows from \eqref{eq:2N5}, once
we write it as a second order wave equation for the primary
wave. Moreover, the data \eqref{eq:ELDat1} take the form \eqref{eq:3},
\begin{equation}
D_k^{(r,s)} = \int_{\Omega_c} {\rm d} \bx \,  \big(\bb^{(r)}(\bx)\big)^T
\bP^{(s)}(t_k, \bx),
\end{equation}
for $t_k = k \tau$, $k = 0,\ldots, 2n-1$ and $r,s = 1, \ldots, m$.

\vspace{0.05in} 
\begin{rem}
Note that although there are two impedances, one for the pressure wave
and one for the shear wave, we have only one reflectivity $q(\bx)$
defined by \eqref{eq:EL15}. This is because once we fix the two
velocities $c_p(\bx)$ and $c_s(\bx)$ we have only one remaining
independent medium dependent parameter, and the two impedances are
related.

This parametrization is consistent with the results in \cite[Section
  2, Fig. 1]{BeylkinBurridge}, which show that when the array is small
i.e., the angle between the outgoing and incoming waves at the array
is small, the leading contribution to the Born approximation comes
from variations of the impedance $\sigma_p(\bx)$ for P-to-P scattering
and $\sigma_s(\bx)$ for S-to-S scattering. Here P stands for pressure
waves and $S$ for shear waves. Mode conversions S-to-P and P-to-S have
a much smaller contribution, proportional to the sine of the angle
between the outgoing and incoming waves, and are neglected in our
approach.
\end{rem}

%=================================================
\section{Numerical results}
\label{sect:Num}
In this section we present numerical results obtained with 
Algorithm \ref{alg:DtB} for the   two dimensional (2D) acoustic wave equation\footnote{Note that the same equation models electromagnetic waves in transverse electric and transverse magnetic modes, with the wave speed $c(\bx)$ and impedance $\sigma(\bx)$ given by \eqref{eq:E5}.} \eqref{eq:S1} and  the 2D  elastic wave  equation \eqref{eq:EL2}. 

\subsection{2D scalar wave problem}
\label{sect:Num1}

\begin{figure}[t!]
\centering
\begin{tabular}{cc}
$\sigma(\bx)$ & $c(\bx)$ \\
\includegraphics[width=0.45\textwidth]{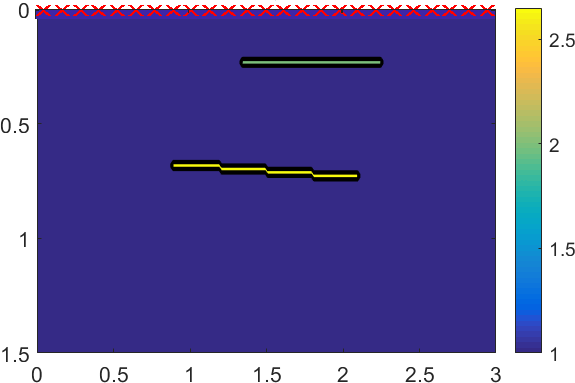} & 
\includegraphics[width=0.41\textwidth]{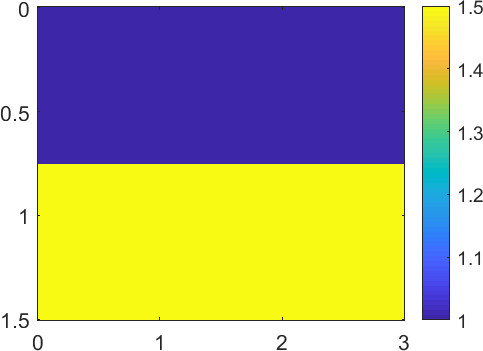}
\end{tabular}
\vspace{-0.12in}\caption{Model acoustic impedance $\sigma(\bx)$ with two inclusions (left) and piecewise-constant velocity
  $c(\bx)$ (right).  The crosses in the left plot indicate the sensor locations.The abscissa and ordinate are in km units. The colorbar displays
  the values of the impedance and velocity normalized by their constant values in the vicinity of the array. }
\label{fig:model_acoustics}
\end{figure}

\begin{figure}[t!]
\centering
\begin{tabular}{cc}
Scattered field & True Born approximant\\
\includegraphics[width=0.45\textwidth]{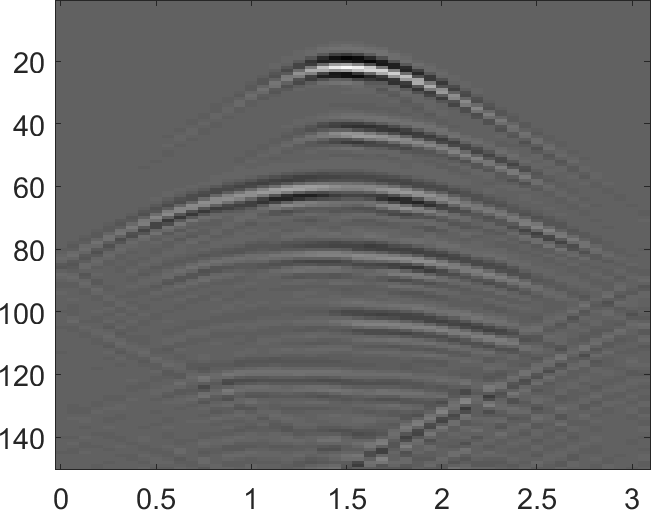} & 
\includegraphics[width=0.45\textwidth]{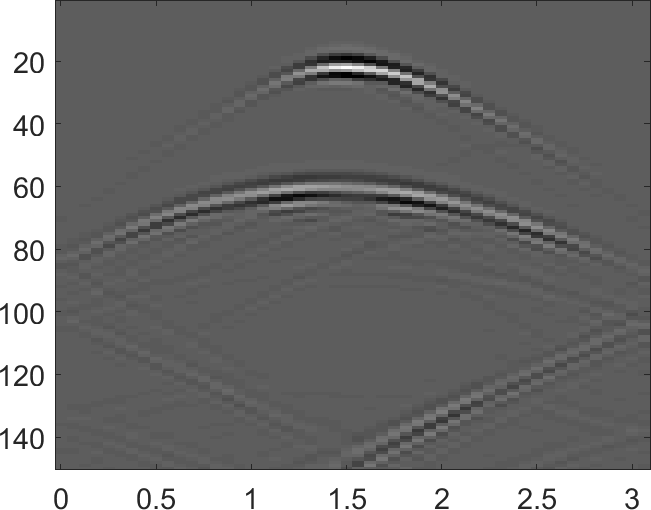}\\
DtB for data with 10\% noise & DtB for noiseless data \\
\includegraphics[width=0.45\textwidth]{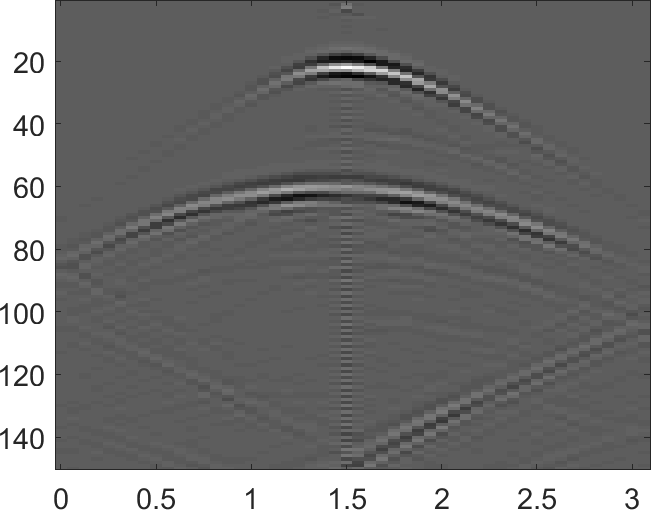} & 
\includegraphics[width=0.45\textwidth]{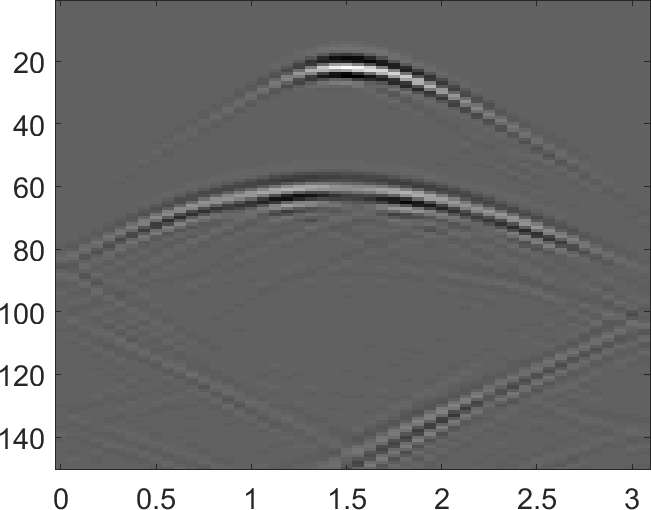}\\
\end{tabular}
\vspace{-0.12in}\caption{The sound wave measured at the array (top left), true Born approximation (top right) and the outputs  of Algorithm \ref{alg:DtB} for data with 10\% noise (bottom left) and noiseless data (bottom right). The abscissa is the sensor location in $km$ units and the ordinate is the index $k$ of the discrete time instants $t_k=k\tau$, with $\tau=0.034s$. The plots are on the same gray scale.}
\label{fig:res_ac}
\end{figure}

Consider the 2D acoustic wave scattering problem in the medium shown in Fig. \ref{fig:model_acoustics}, where we plot
the  acoustic impedance $\sigma(\bx)$ and wave speed $c(\bx)$ normalized by their constant values at the array. The data $\bD$ and its Born approximation $\bD^{{\rm Born}}$ are obtained using finite-difference time-domain simulations with time step close to the CFL limit. The separation between the $m_a=50$ sensors and the time sampling rate $\tau$ are close to the Nyquist limit. 
The computational domain is the rectange shown in Fig. \ref{fig:model_acoustics}, with sides in km units. At the top boundary 
we have homogeneous Neumann conditions and at the remaining part of the boundary we have homogenenous Dirichlet conditions.

In the top left plot of Fig. \ref{fig:res_ac} we display the raw scattered data due to the excitation from the sensor located  in the middle of the array. The single scattering pattern (primaries)  consist of two reflections from the two thin inclusions and the reflections from the domain boundary, as seen from the Born approximation displayed in the top right plot in Fig. \ref{fig:res_ac}. The other reflections in the top left plot are caused by  multiple scattering  between the inclusions and/or the boundary. The output of the Algorithm \ref{alg:DtB} is displayed in the bottom row of Fig. \ref{fig:res_ac}, for noiseless data in the right plot and data contaminated with 10\% additive  i.i.d.  Gaussian noise 
in the left plot. Both results are almost the same as the true Born approximant. 

%=================================================
\subsection{2D elastic wave problem}
\label{sect:Num2}

The data $\bD$ and the Born approximation $\bD^{{\rm Born}}$ in this section are obtained by solving the 
elastic wave equation \eqref{eq:EL4} with boundary conditions 
\eqref{eq:EL5N}--\eqref{eq:EL5} and initial conditions \eqref{eq:EL6} using  a finite-difference time-domain method with time step satisfying the CFL condition for the shear wave.  The spacing between the $m_a=25$ sensors in the array  and the time sampling rate are close to the Nyquist limit for the shear wave. Therefore, the pressure wave is spatially oversampled. At each sensor location we consider both orientations of the external force source and measure both the horizontal and vertical  velocities. 

\begin{figure}[t!]
\centering
\includegraphics[width=0.45\textwidth]{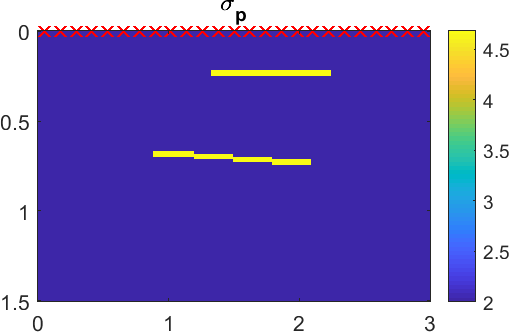}
\vspace{-0.12in}\caption{The pressure wave impedance $\sigma_p(\bx)$ modeling two thin inclusions. The axes are in km units. Crosses indicate the sensor locations. The colorbar shows the impedance normalized by its value at the array. }
\label{fig:elas_model_fr}
\end{figure}

\begin{figure}[t!]
\centering
\begin{tabular}{ccc}
Scattered field & True Born approximant& DtB output\\
\hspace{-0.12in} \includegraphics[width=0.327\textwidth]{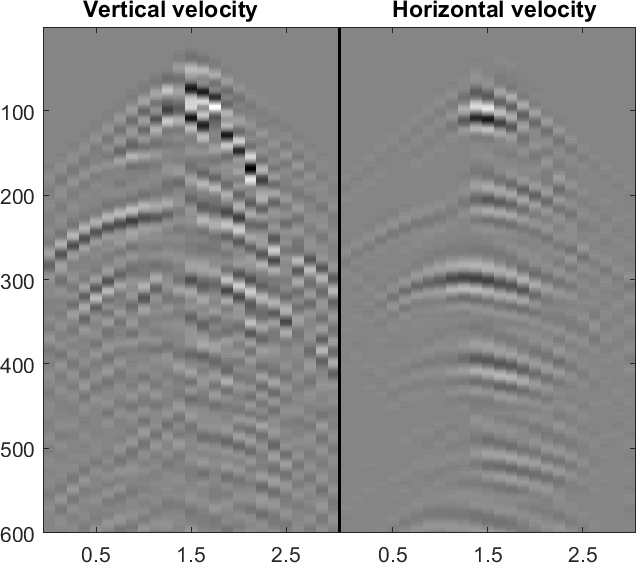} & \hspace{-0.18in}
\includegraphics[width=0.327\textwidth]{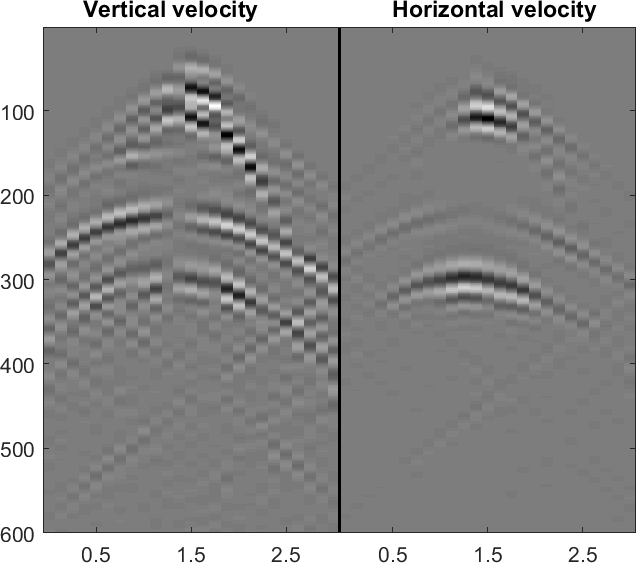}& \hspace{-0.18in}
\includegraphics[width=0.327\textwidth]{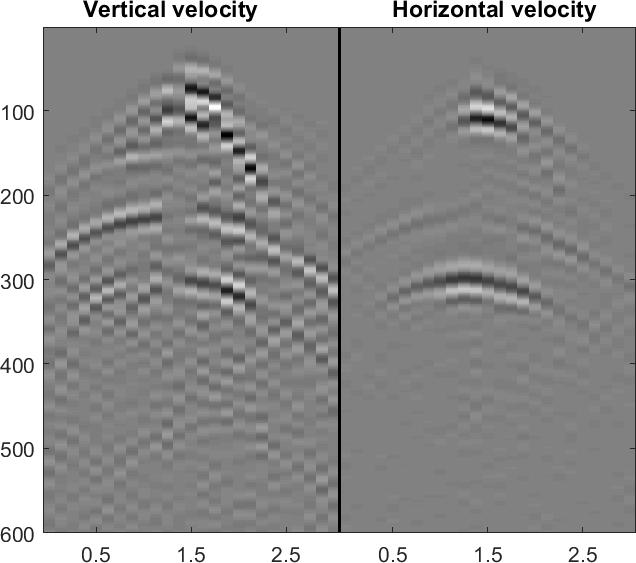}
\end{tabular}
\vspace{-0.12in}\caption{The raw scattered field (left), the true Born approximant (middle) and the output of Algorithm \ref{alg:DtB} (right).  The abscissa is the location of sensors in km units. These measure both the horizontal and vertical velocity components.  The ordinate is the index $k$ of the discrete time instants $t_k=k\tau$, with $\tau=0.034s$. Because the horizontal response is stronger, we amplified the vertical response by the factor 8. This amplification highlights numerical artifacts  in the vertical velocity plots. The plots are on the same gray scale.}
\label{fig:elas_res_fr}
\end{figure}

In the first  simulation we consider  two thin inclusions embedded in a homogeneous background. The wave speeds $c_p$ and 
$c_s$ are constant, satisfying $c_p = 2 c_s$, so the shear wave impedance is $\sigma_s(\bx) = \sigma_p(\bx)/2$. 
In Fig. \ref{fig:elas_model_fr} we display the normalized (by the value of $\sigma_s(\bx)$ at the array) pressure wave 
impedance $\sigma_p(\bx)$.

We present results  for the horizontal force exerted by the source located in the middle of the array.
The left plot in Fig. \ref{fig:elas_res_fr} displays the raw array data, which contain the primary arrivals and the multiply 
scattered waves between the inclusions and/or the boundary of the domain.
Since there are two waves traveling with different speeds, the primaries consist of two pressure waves and two shear 
waves reflected from the two inclusions and the domain boundaries. These can be seen  in the middle plot in 
Fig. \ref{fig:elas_res_fr}. Due to the excitation and the nearly layered medium, the dominant response is the 
horizontal velocity of the shear waves. The weaker vertical velocity response is amplified in the plots by the factor $8$,
in order to display it on the same gray scale as the horizontal response.  The output of the Algorithm \ref{alg:DtB} 
is shown in the right plot of Fig. \ref{fig:elas_res_fr}. It is basically the same as the Born approximation, aside from 
some artifacts in the vertical velocity plot. 

\begin{figure}[t!]
\centering
\includegraphics[width=0.45\textwidth]{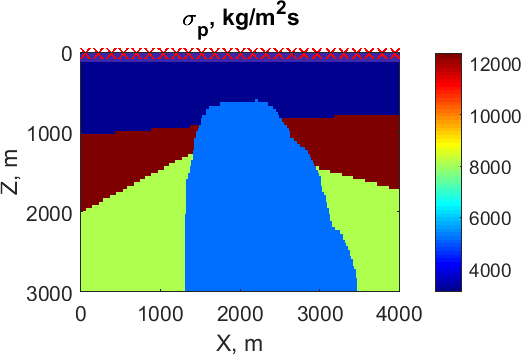}
\vspace{-0.12in}\caption{Salt dome impedance model of the pressure wave impedance $\sigma_p(\bx)$. The abscissa and 
ordinate are in m units. The colorbar is in ${\rm kg}/({m}^2s)$ units. Crosses indicate the sensor locations.}
\label{fig:elas_model_sd}
\end{figure}

\begin{figure}[t!]
\centering
\begin{tabular}{ccc}
Scattered field & True Born approximant& DtB output\\
\hspace{-0.12in}
\includegraphics[width=0.327\textwidth]{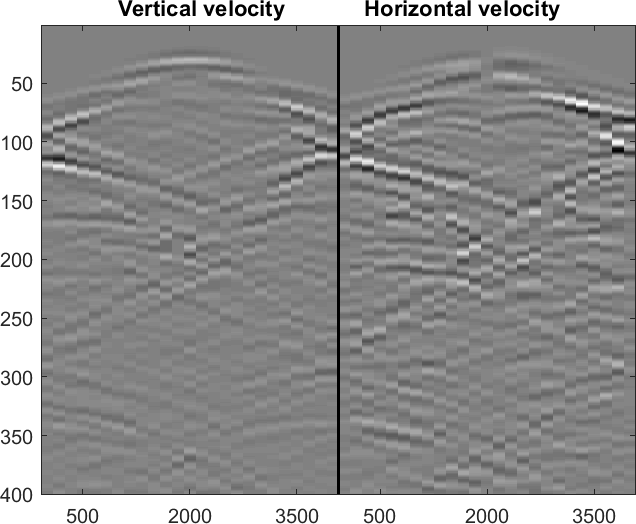} & \hspace{-0.18in}
\includegraphics[width=0.327\textwidth]{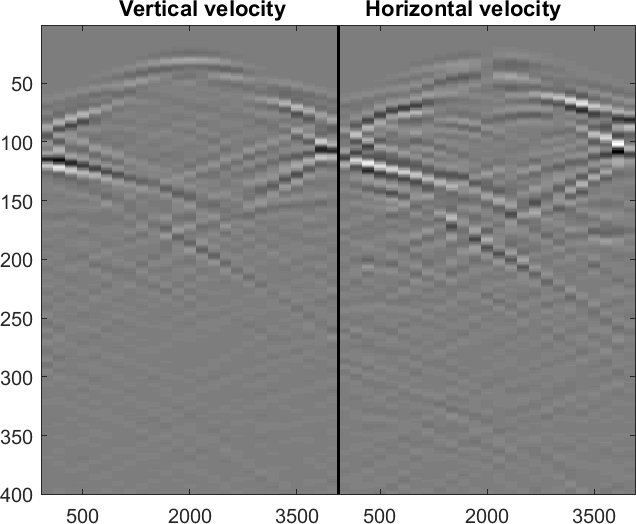}&  \hspace{-0.18in}
\includegraphics[width=0.327\textwidth]{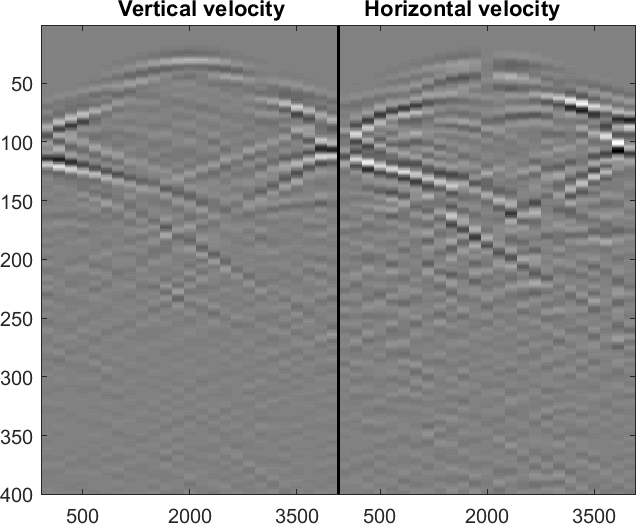}
\end{tabular}
\vspace{-0.12in}\caption{Raw scattered field (left), the true Born approximant (middle) and the output of DtB transform (right). The abscissa is location of sensors measuring horizontal and vertical velocity components, and the ordinate is the index $k$ of the discrete time instants $t_k=k\tau$, $\tau=0.034s$. The plots are on the same gray scale.}
\label{fig:elas_res_sd}
\end{figure}
}

The second simulation is motivated by the application of seismic exploration, and models a salt dome in a layered formation, 
as shown in Fig. \ref{fig:elas_model_sd}. The pressure and shear wave speeds are constant and equal to 3100m/s and 1800m/s, respectively. The background pressure wave impedance is homogeneous and equal to  3100${\rm kg}/({\rm m}^2{\rm s})$. 
As in the first simulation, we consider both horizontal and vertical external forces and the data consist of both the horizontal 
and vertical velocities. In Fig. \ref{fig:elas_res_sd} we display the raw and processed data for the vertical force exerted from 
the source in the center of the array. In the left plot we show the raw data. Since the medium is far from layered, both responses 
are of the same order and there is no amplification factor in the plots. The primary reflections can be seen in the middle plot and 
the  output of the Algorithm \ref{alg:DtB} is shown in the right plot. It matches the true Born approximation displayed in the 
middle plot in Fig. \ref{fig:elas_res_sd}. 

%=================================================
\section{Summary}
\label{sect:sum}
We introduced a robust algorithm for  nonlinear processing of data gathered by an active array of sensors, 
which seeks to determine an unknown medium by probing it with pulses and measuring the resulting waves. 
These waves depend nonlinearly on the variations in the medium, modeled by an unknown  reflectivity function. 
Many imaging methodologies ignore this nonlinearity and operate under the linear, single scattering (Born) assumption.  
This is adequate for a weak reflectivity. However, in strongly scattering media the nonlinear (multiple scattering) effects 
are significant and images based on the Born approximation have unwanted artifacts. The  algorithm introduced in this 
paper seeks to map the array data gathered for an arbitrary (possibly large) reflectivity  to the single scattering (Born) 
data, which can then be used by any linear imaging method. This mapping is called the Data-to-Born (DtB) transformation.

The algorithm is based on a data driven reduced order model which consists of a proxy wave propagator operator. The true 
wave propagator  maps the wave from a given state at time $t$ to a future state at time $t + \tau$, where $\tau$ is the time 
sampling of the measurements. The proxy wave propagator is a projection of the true wave propagator on the space 
spanned by the snapshots of the wave field at the discrete times $t_k = k \tau$, for $k = 0, \ldots, n-1$. It is constructed 
directly from the array data at $t_k = k \tau$, for $k = 0, \ldots, 2n-1$,  with no knowledge of the reflectivity in the medium. 

Our definition of the reflectivity function is based on the known fact that  the main contribution to the Born approximation 
for a small array is due to variations of the logarithm of the impedance of the medium.  The wave speed, which determines 
the kinematics (travel times) of the wave is assumed known. The DtB transformation is obtained using a factorization of the 
proxy wave propagator in two operators that have an approximately affine dependence on the reflectivity function. 
This allows the computation of the Fr\'{e}chet derivative of the reflectivity to data mapping, which defines the Born 
approximation of the measurements. 

The algorithm is developed for a generic hyperbolic system and we showed how it applies to the three types of linear 
waves: sound, electromagnetic and elastic.  Because it consists of a sequence of algebraic operations that can be 
performed without knowing the exact wave propagation model, the algorithm  is versatile and can be used as a black-box 
tool for any of these waves. 

To ensure robustness for noisy data, we identified the unstable step in the algorithm and introduced a regularization 
procedure based on a spectral truncation of the data driven reduced order model. This regularization balances 
numerical stability and accuracy of data fitting,  up to the order of the standard deviation of the noise. 
The performance of the algorithm is assessed with numerical simulations for both sound and elastic waves in 
two dimensions.

%=================================================
\section*{Acknowledgements}
This material is based upon research supported in part by the
U.S. Office of Naval Research under award number N00014-17-1-2057 to
Borcea and Mamonov. Borcea also acknowledges support from the AFOSR
award FA9550-18-1-0131 and 
Mamonov acknowledges support from the National Science Foundation Grant DMS-1619821.
The research was done in  part at the Institute for Computational 
and Experimental Research in Mathematics in Providence, RI, supported by the National Science Foundation under 
Grant No. DMS-1439786, where the authors were in residence during the Fall 2017 semester.
We thank Dr. Smaine Zeroug from Schlumberger-Doll Research for his suggestions on elasticity benchmarks.

%========================================================================
\bibliography{biblio} 
\bibliographystyle{siam}

\end{document}

%% file: mathmacros.tex
\DeclareMathAlphabet{\itbf}{OML}{cmm}{b}{it}
 \DeclareMathAlphabet\mathbfcal{OMS}{cmsy}{b}{n}

\renewcommand{\hat}{\widehat}
\renewcommand{\tilde}{\widetilde}
\def\RR{\mathbb{R}}
\def\bx{{{\itbf x}}}
\def\vx{{\vec{\itbf x}}}
\def\bu{{{\itbf u}}}
\def\bb{{\itbf b}}

\def\bP{{\itbf P}}

\def\bv{{\itbf v}}

\def\dP{{d_{_{\hspace{-0.01in}P}}}}

\def\bV{{\itbf V}}
\def\bR{{\itbf R}}

\def\bE{{\itbf E}}

\def\bL{{\bf L}}
\def\bM{{\bf M}}
\def\bSig{\boldsymbol{\mathfrak S}}
\def\cbS{\boldsymbol{\mathcal S}}

\def\bS{{\bf S}}

\def\bZ{\boldsymbol{\mathcal{Z}}}

\def\bQQ{{\bf Q}}
\def\bK{\boldsymbol{\mathcal K}}
\def\bD{{\bf D}}
\def\cP{{\mathscr P}}
\def\cbP{{\boldsymbol{\cP}}_{\hspace{-0.03in}q}}
\def\cbPn{\cbP^{\hspace{-0.01in}(n)}}
\def\cbPno{{\boldsymbol{\cP}}_0^{\hspace{-0.01in}(n)}}
\def\bbn{\bb^{(n)}}
\def\cbPN{\cbP^{\hspace{-0.01in}(N)}}
\def\cbLN{{\boldsymbol{\cal L}}^{\hspace{-0.01in}(N)}}
\def\cbLNT{{\boldsymbol{\cal L}}^{\hspace{-0.01in}(N)T}}
\def\cbLn{{\boldsymbol{\cal L}}^{\hspace{-0.01in}(n)}}
\def\tbLn{{\widetilde{\boldsymbol{\cal L}}}^{\hspace{-0.01in}(n)}}
\def\tbLnT{{\widetilde{\boldsymbol{\cal L}}}^{\hspace{-0.01in}(n)T}}
\def\cbLnT{{\boldsymbol{\cal L}}^{\hspace{-0.01in}(n)T}}
\def\bbN{\bb^{(N)}}
\def\cT{{\mathcal T}}
\def\ep{\epsilon}
\def\vep{\varepsilon}

\def\bet{{\boldsymbol \eta}}
\def\bzet{{\boldsymbol \zeta}}

\def\bpi{{\boldsymbol \pi}}
\def\bett{\widetilde{\bet}}
\def\bGamma{{\boldsymbol{\Gamma}}}

\def\thP{\widetilde {\widehat
   \bP}^{\mbox{\raisebox{-0.1in}{\scriptsize $({n})$}}}}
\def\cPq{{\cP_{\hspace{-0.03in}q}}}
\def\bSt{\boldsymbol{\mathfrak{T}}}
\def\St{{\mathfrak{T}}}
\def\bcD{\boldsymbol{\mathcal D}}